\documentclass[12pt]{article}
\usepackage{amssymb,amsthm,amsmath,amsfonts,latexsym,tikz,hyperref,skak}
\usepackage[notref, notcite]{showkeys}
\usepackage[hmargin=.8in,vmargin=.9in]{geometry}

\newtheorem{thm}{Theorem}[section]
\newtheorem{prop}[thm]{Proposition}
\newtheorem{cor}[thm]{Corollary}
\newtheorem{lem}[thm]{Lemma}
\newtheorem{conj}[thm]{Conjecture}
\newtheorem{exa}[thm]{Example}

\newcommand{\cQ}{{\cal Q}}
\newcommand{\da}{\hs{-2pt}\downarrow}
\newcommand{\lf}{\lfloor}
\newcommand{\rf}{\rfloor}

\newcommand{\omh}{\hat{\om}}
\newcommand{\omc}{\check{\om}}

\DeclareMathOperator{\coinv}{coinv}
\DeclareMathOperator{\area}{area}
\DeclareMathOperator{\dinv}{dinv}
\DeclareMathOperator{\bounce}{bounce}

\newcommand{\ben}{\begin{enumerate}}
\newcommand{\een}{\end{enumerate}}
\newcommand{\ble}{\begin{lem}}
\newcommand{\ele}{\end{lem}}
\newcommand{\bth}{\begin{thm}}
\renewcommand{\eth}{\end{thm}}
\newcommand{\bpr}{\begin{prop}}
\newcommand{\epr}{\end{prop}}
\newcommand{\bco}{\begin{cor}}
\newcommand{\eco}{\end{cor}}
\newcommand{\bcon}{\begin{conj}}
\newcommand{\econ}{\end{conj}}
\newcommand{\bde}{\begin{defn}}
\newcommand{\ede}{\end{defn}}
\newcommand{\bex}{\begin{exa}}
\newcommand{\eex}{\end{exa}}
\newcommand{\barr}{\begin{array}}
\newcommand{\earr}{\end{array}}
\newcommand{\btab}{\begin{tabular}}
\newcommand{\etab}{\end{tabular}}
\newcommand{\beq}{\begin{equation}}
\newcommand{\eeq}{\end{equation}}
\newcommand{\bea}{\begin{eqnarray*}}
\newcommand{\eea}{\end{eqnarray*}}
\newcommand{\bal}{\begin{align*}}
\newcommand{\bce}{\begin{center}}
\newcommand{\ece}{\end{center}}
\newcommand{\bpi}{\begin{picture}}
\newcommand{\epi}{\end{picture}}
\newcommand{\bpp}{\begin{picture}}
\newcommand{\epp}{\end{picture}}
\newcommand{\bfi}{\begin{figure} \begin{center}}
\newcommand{\efi}{\end{center} \end{figure}}
\newcommand{\bprf}{\begin{proof}}
\newcommand{\eprf}{\end{proof}\medskip}
\newcommand{\capt}{\caption}
\newcommand{\bsl}{\begin{slide}{}}
\newcommand{\esl}{\end{slide}}
\newcommand{\bfr}{\begin{frame}}
\newcommand{\efr}{\end{frame}}

\newcommand{\hqed}{\hfill \qed}

\newcommand{\eqqed}[1]{$\rule{1ex}{0ex}\hfill{\dil#1}\hfill\qed$}

\newcommand{\ol}{\overline}

\newcommand{\hs}[1]{\hspace{#1}}
\newcommand{\hso}[1]{\hspace{-1pt}}
\newcommand{\vs}[1]{\vspace{#1}}
\newcommand{\qmq}[1]{\quad\mbox{#1}\quad}

\newcommand{\emp}{\emptyset}

\newcommand{\sbs}{\subset}
\newcommand{\sbe}{\subseteq}

\newcommand{\spe}{\supseteq}
\newcommand{\setm}{\setminus}

\newcommand{\gauss}[2]{\genfrac{[}{]}{0pt}{}{#1}{#2}}

\def\<{\langle}
\def\>{\rangle}

\newcommand{\ree}[1]{(\ref{#1})}

\newcommand{\ra}{\rightarrow}

\newcommand{\al}{\alpha}
\newcommand{\be}{\beta}

\newcommand{\ep}{\epsilon}

\newcommand{\om}{\omega}

\newcommand{\ze}{\zeta}

\newcommand{\De}{\Delta}

\newcommand{\bbZ}{{\mathbb Z}}

\newcommand{\cB}{{\cal B}}

\newcommand{\cR}{{\cal R}}

\newcommand{\fS}{{\mathfrak S}}

\newcommand{\Bb}{\ol{B}}

\newcommand{\pib}{\ol{\pi}}

\DeclareMathOperator{\inv}{inv}

\DeclareMathOperator{\wt}{wt}

\newcommand{\dil}{\displaystyle}

\begin{document}
\pagestyle{plain}

\title{$m$-level rook placements
}
\author{Kenneth Barrese\\[-5pt]
\small Department of Mathematics, Michigan State University,\\[-5pt]
\small East Lansing, MI 48824-1027  {\tt baressek@math.msu.edu}\\[5pt]
Nicholas Loehr \thanks{This work was partially supported by a grant from the Simons Foundation
   (\#244398 to Nicholas Loehr).}\\[-5pt]
\small Department of Mathematics, Virginia Tech\\[-5pt]
\small Blacksburg, VA 24061-0123 {\tt nloehr@vt.edu}\\[-5pt]
\small and\\[-5pt]
\small Department of Mathematics, United States Naval Academy\\[-5pt]
\small Annapolis, MD 21402-5002 {\tt loehr@usna.edu}\\[5pt]
Jeffrey Remmel\\[-5pt]
\small Department of Mathematics, UCSD\\[-5pt]
\small La Jolla, CA, 92093-0112 {\tt jremmel@ucsd.edu}\\[5pt]
Bruce E. Sagan\\[-5pt]
\small Department of Mathematics, Michigan State University,\\[-5pt]
\small East Lansing, MI 48824-1027 {\tt sagan@math.msu.edu}
}

\date{\today\\[10pt]
	\begin{flushleft}
	                                       \\[5pt]
	\end{flushleft}}

\maketitle

\begin{abstract}
Goldman, Joichi, and White proved a beautiful theorem showing that the falling factorial generating function for the rook numbers of a Ferrers board factors over the integers.  Briggs and Remmel studied an analogue of rook placements where rows are replaced by sets of $m$ rows called levels.  They proved a version of the factorization theorem in that setting, but only for certain Ferrers boards.  We generalize this result to any Ferrers board as well as giving a $p,q$-analogue.  We also consider a dual situation involving weighted file placements which permit more than one rook in the same row.  In both settings, we discuss properties of the resulting equivalence classes such as the number of elements in a class.  In addition, we prove analogues of a theorem of Foata and Sch\"utzenberger giving a distinguished representative in each class as well as make connections with the $q,t$-Catalan numbers.  We end with some  open questions raised by this work.
\end{abstract}

%%%%%%%%%%%%%%%%%%%%%%%%%%%%%%%%
%
% 	INTRODUCTION 	
%
%%%%%%%%%%%%%%%%%%%%%%%%%%%%%%%%

\section{Introduction}
\label{i}

\bfi
\begin{tikzpicture}
\draw(-1,2) node {$Q=$};
\draw(2,5) node{$\vdots$};
\draw(5,2) node{$\dots$};
\draw(0,0) grid(4,4);
\fill(1.5,2.5) circle(.1);
\end{tikzpicture}
\hs{50pt}
\begin{tikzpicture}
\foreach \y in {0,1} 
   \draw (0,\y)--(3,\y);
\foreach \y in {2,3,} 
   \draw (1,\y)--(3,\y);
\foreach \x in {1,2,3}
   \draw (\x,0)--(\x,3);
\draw (0,0)--(0,1);
\draw(-2,1.5) node{$B=(1,3,3)=$};
\end{tikzpicture}
\capt{The quadrant and a board \label{qb}}
\efi

Our point of departure will be the famous Factorization Theorem of Goldman, Joichi, and White~\cite{gjw:rtI}.  To state it, we first need to set up some standard notation for rook theory.  Consider the tiling, $Q$, of the first quadrant with unit squares (also called ``cells" or ``boxes") displayed on the left in Figure~\ref{qb}.  We give each square the coordinates of its northeast corner, so the circle in the diagram of $Q$ is in cell $(2,3)$.  A \emph{board} is  finite subset $B\sbs Q$.  We will be particularly interested in boards associated with (integer) partitions.  A \emph{partition} is a weakly increasing sequence $(b_1,\ldots,b_n)$ of nonnegative integers.  We will use the same notation for the corresponding  \emph{Ferrers board} $B=(b_1,\ldots,b_n)$ which consists of  the $b_j$ lowest squares in column $j$ for $1\le j\le n$.  The board $B=(1,3,3)$ is shown on the right in Figure~\ref{qb} and it will be used as our running example for this section.

For any board, $B$, a \emph{rook placement} is a subset $P\sbe B$ such that no two squares of $P$ are in the same row or column.  The elements of $P$ are usually called \emph{rooks}.  We let 
$$
r_k(B)=\text{the number of rook placements $P\sbe B$ with $k$ rooks.}
$$
Note that we always have $r_0(B)=1$ and $r_1(B)=|B|$ where $|\cdot|$ denotes cardinality.  For $B=(1,3,3)$ we have $r_0(B)=1$, $r_1(B)=7$, $r_2(B)=10$, $r_3(B)=2$, and $r_k(B)=0$ for $k\ge4$.  We wish to consider the generating function for these integers in the variable $x$ and using the basis of \emph{falling factorials}
$$
x\da_n = x(x-1)(x-2)\dots(x-n+1)
$$
for $n\ge0$.  This brings us to the Factorization Theorem.
\bth[Factorization Theorem~\cite{gjw:rtI}]
\label{gjw}
If  $B=(b_1,\dots,b_n)$ is any Ferrers board then

\vs{10pt}

\eqqed{
\sum_{k=0}^n r_k(B) x\da_{n-k} = \prod_{j=1}^n (x+b_j-j+1).
}
\eth

Motivated by connections to wreath products of cyclic groups with symmetric groups, Briggs and Remmel~\cite{br:mrn} considered rook placements where rows are replaced by sets of rows called levels.  Fix a positive integer $m$.  Partition the rows of $Q$ into levels where the \emph{$i$th level} consists of rows $(i-1)m+1,(i-1)m+2,\dots,im$.  The situation for $m=2$ is shown on the left in Figure~\ref{lrp} where the boundaries between the levels have been thickened.  

Given a board, $B$, an \emph{$m$-level rook placement} (called an \emph{$m$-rook placement} by Briggs and Remmel) is $P\sbe B$ where no two elements of $R$ are in the same level or the same column.  Note that when $m=1$ we recover the ordinary notion of a rook placement.  By way of example, in Figure~\ref{lrp}, 
the  placement  on the middle board is a $2$-level rook placement while the one on the right is not since it has two rooks in the first level.
We let
$$
r_{k,m}(B)=\text{the number of $m$-level rook placements $P\sbe B$ with $k$ rooks.}
$$
In general, we will add a subscript $m$ to quantities when considering their $m$-level equivalents.    For $B=(1,3,3)$ we have 
$r_{0,2}(B)=1$, $r_{1,2}(B)=7$, $r_{2,2}(B)=6$, and $r_k(B)=0$ for $k\ge3$. 

\bfi
\begin{tikzpicture}
\draw(-1.5,1) node{level $1$};
\draw(-.5,1) node{\scalebox{4}{$\{$}};
\draw(-1.5,3) node{level $2$};
\draw(-.5,3) node{\scalebox{4}{$\{$}};
\draw(2,5) node{$\vdots$};
\draw(5,2) node{$\dots$};
\draw(0,0) grid(4,4);
\draw[very thick] (0,0)--(4,0) (0,2)--(4,2) (0,4)--(4,4);
\end{tikzpicture}
\hs{40pt}
\begin{tikzpicture}
\foreach \y in {0,1} 
   \draw (0,\y)--(3,\y);
\foreach \y in {2,3,} 
   \draw (1,\y)--(3,\y);
\foreach \x in {1,2,3}
   \draw (\x,0)--(\x,3);
\draw (0,0)--(0,1);
\draw[very thick] (1,2)--(3,2) (0,0)--(3,0);
\draw(.5,.5) node {\scalebox{2}{\symrook}};
\draw(2.5,2.5) node {\scalebox{2}{\symrook}};
\end{tikzpicture}
\hs{40pt}
\begin{tikzpicture}
\foreach \y in {0,1} 
   \draw (0,\y)--(3,\y);
\foreach \y in {2,3,} 
   \draw (1,\y)--(3,\y);
\foreach \x in {1,2,3}
   \draw (\x,0)--(\x,3);
\draw (0,0)--(0,1);
\draw[very thick] (1,2)--(3,2) (0,0)--(3,0);
\draw(.5,.5) node {\scalebox{2}{\symrook}};
\draw(2.5,1.5) node {\scalebox{2}{\symrook}};
\end{tikzpicture}
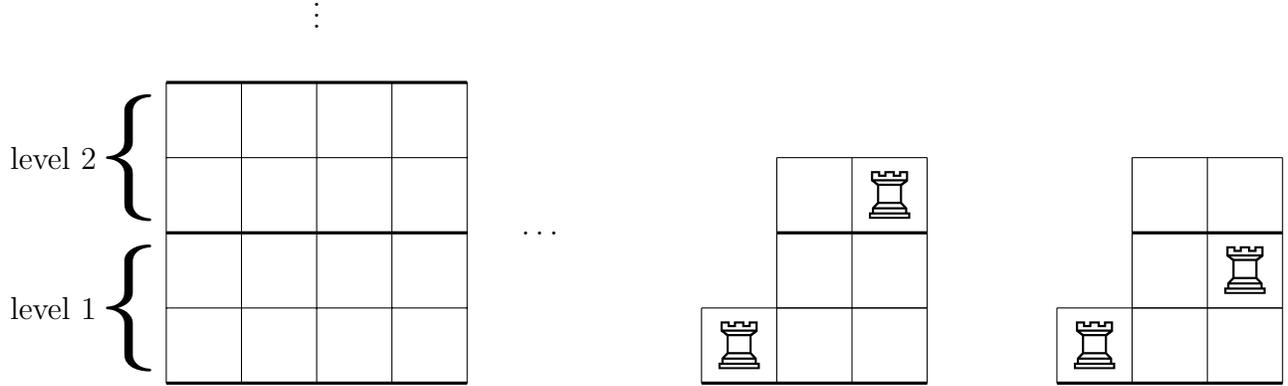
\capt{Levels and rook placements \label{lrp}}
\efi

To state the Briggs-Remmel generalization of Theorem~\ref{gjw}, we need a few  more concepts.  One is of an \emph{$m$-falling factorial} which is
$$
x\da_{n,m}=x(x-m)(x-2m)\dots(x-(n-1)m).
$$
Another is the \emph{$m$-floor function} defined by
$$
\lf n \rf_m = \text{the largest multiple of $m$ less than or equal to $n$}
$$
for any integer $n$.  As an example $\lf 17 \rf_3 = 15$ since $15\le 17 < 18$.  Finally, define a \emph{singleton board} to be a Ferrers board $B=(b_1,\dots,b_n)$ such that $\lf b_i \rf_m \neq b_i$ implies $\lf b_{i+1} \rf_m > \lf b_i \rf_m$ for all $i<n$.  These were called \emph{$m$-Ferrers boards} in~\cite{br:mrn}.  The reason for our terminology will be explained when we introduce the concept  of a zone.

\bth[\cite{br:mrn}]
\label{br}
If $B=(b_1,\dots,b_n)$ is a singleton board then

\vs{10pt}

\eqqed{
\sum_{k=0}^n r_{k,m}(B) x\da_{n-k,m} = \prod_{j=1}^n (x+b_j-(j-1)m).
}
\eth

Our first goal is to remove the singleton board restriction and prove a generalization of this theorem for any Ferrers board.  This will be done in the next section.  In Section~\ref{pqa}, we give a $p,q$-analogue of our result using statistics related to the inversion number of a permutation.  Call boards $B, B'$  \emph{$m$-level rook equivalent} if $r_{k,m}(B)=r_{k,m}(B')$ for all $k$.  In Section~\ref{re} we extend to all $m$ a theorem of Foata and Sch\"utzenberger~\cite{fs:rpf} giving a distinguished member of each $1$-level rook equivalence class.   Goldman, Joichi and White used the Factorization Theorem to enumerate the number of Ferrers boards $1$-level rook equivalent to a given board.  In Section~\ref{esb} we generalize this formula to count $m$-level  rook equivalent singleton boards for arbitrary $m$. 
The rest of the paper is devoted to the following dual problem.  Rather than changing the product side of Theorem~\ref{br}, keep the same product for all Ferrers boards and expand it in the $m$-falling factorial basis.  What do the coefficients count?  We show in Section~\ref{fp} that they are generating functions for certain weighted file placements, where such placements allow more than one rook in a given row.  The next two sections investigate properties of the corresponding equivalence classes.  In particular, in Section~\ref{wecs} we count the number of boards in a given class and show how this formula can be obtained using ideas from the theory of $q,t$-Catalan numbers.  The last section contains some open questions related to our work.

%%%%%%%%%%%%%%%%%%%%%%%%%%%%%%%%
%
% 	THE $m$-FACTORIZATION THEOREM 	
%
%%%%%%%%%%%%%%%%%%%%%%%%%%%%%%%%

\section{The $m$-Factorization Theorem}
\label{sec:mft}

\bfi
\begin{tikzpicture}
\fill[lightgray](0,0) rectangle (3,1);
\fill[lightgray](2,1) rectangle (3,2);
\fill[lightgray](4,3) rectangle (5,5);
\fill[lightgray](5,6) rectangle (6,7);
\draw (0,0) -- (6,0) (0,1)--(6,1) (2,2)--(6,2) (3,3)--(6,3) (4,4)--(6,4) (4,5)--(6,5)  (5,7)--(6,7)
(0,0)--(0,1) (1,0)--(1,1) (2,0)--(2,2) (3,0)--(3,3) (4,0)--(4,5) (5,0)--(5,7) (6,0)--(6,7)
;
\draw[very thick] (0,0)-- (6,0) (3,3)--(6,3) (5,6)--(6,6) (3,0)--(3,3) (5,0)--(5,7);
\end{tikzpicture}
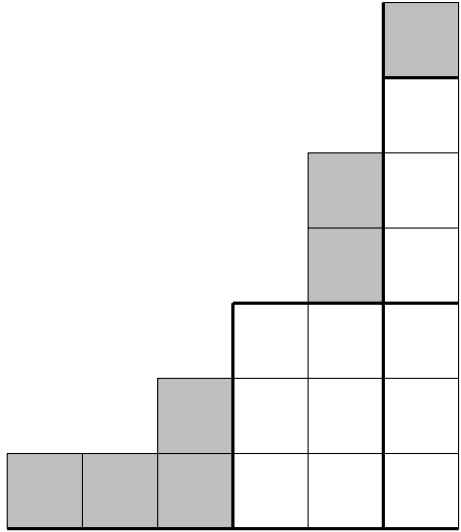
\capt{The zones of $(1,1,2,3,5,7)$ when $m=3$ \label{z}}
\efi

In order to generalize Theorem~\ref{br} to all Ferrers boards, it will be convenient to break a board up into zones depending on the lengths of the columns.  Given integers $s,t$, the interval from $s$ to $t$ will be denoted $[s,t]=\{s,s+1,\dots,t\}$.  An \emph{$m$-zone}, $z$, of a board $B=(b_1,\ldots,b_n)$ is  a maximal interval $[s,t]$ such that 
$\lf b_s \rf_m =\lf b_{s+1} \rf_m = \dots =\lf b_t \rf_m$.
To illustrate this concept, consider $m=3$ and the board $B=(1,1,2,3,5,7)$ shown in Figure~\ref{z}.  In this case the zones are $z_1=[1,3]$ since $\lf b_1 \rf_3=\lf b_2 \rf_3 =\lf b_3 \rf_3=0$, 
$z_2=[4,5]$ since $\lf b_4 \rf_3=\lf b_5 \rf_3=3$, and $z_3=[6,6]$ since $\lf b_6 \rf_3 = 6$.  The zones in Figure~\ref{z} are separated by thick lines (as are the levels).  Note that a Ferrers board is a singleton board if and only if each zone contains at most one column whose length is not a multiple of $m$.  This is the reason for our choice of terminology.

In addition to taking $m$-floors, we will have to consider remainders modulo $m$.  Given an integer $n$, we denote its remainder on division by $m$ by $\rho_m(n)=n-\lf n\rf_m$.  If $z$ is a  zone of a Ferrers board $B=(b_1,\dots,b_n)$ then its \emph{$m$-remainder} is 
$$
\rho_m(z) = \sum_{j\in z} \rho_m(b_j).
$$
In Figure~\ref{z}, the boxes corresponding to the $3$-remainders of the zones are shaded.  In particular $\rho_3(z_1)=1+1+2=4$, $\rho_3(z_2)=0+2=2$, and $\rho_3(z_3)=1$.  We are now in a position to state and prove our generalization of Theorem~\ref{br}.
\bth[$m$-Factorization Theorem]
\label{mft}
If $B=(b_1,\dots,b_n)$ is any Ferrers board then
$$
\sum_{k=0}^n r_{k,m}(B) x\da_{n-k,m} = \prod_{j=1}^n 
\begin{cases}
x+\lf b_j \rf_m -(j-1)m + \rho_m(z) &\text{if $j$ is the last index in its zone $z$,}\\
x+\lf b_j \rf_m -(j-1)m &\text{otherwise.}
\end{cases}
$$
\eth
\bprf
Since this is a polynomial identity, it suffices to prove it for an infinite number of values for $x$.  We will do so when $x$ is a nonnegative multiple of $m$.  Consider the board $B_x$ derived from $B$ by adding an $x\times n$ rectangle below $B$.  Figure~\ref{B_x} shows a schematic representation of $B_x$.  Note that since $x$ is a multiple of $m$, the zones and remainders of $B$ and $B_x$ are the same.   We will show that both the sum and the product count the number of $m$-rook placements on $B_x$ consisting of $n$ rooks.

\bfi
\begin{tikzpicture}
\draw (0,0)--(4,0)--(4,8)--(2,8)--(2,6)--(0,6)--(0,0) (0,4)--(4,4)
;
\draw[very thick, <->] (0,-.5)--(4,-.5);
\draw[very thick, <->] (-.5,0)--(-.5,4);
\draw(2,-1) node {$n$};
\draw(-1,2) node {$x$};
\draw(2.5,5.5) node {$B$};
\draw(-1.5,4) node {$B_x=$};
\end{tikzpicture}
\capt{ The board $B_x$ \label{B_x}}
\efi

For the sum side, note that any placement of $n$ rooks on $B_x$ must have $k$ rooks in $B$ and $n-k$ rooks in the rectangle for some $0\le k\le n$.  By definition, $r_{k,m}(B)$ counts the number of placements on $B$.  Once these rooks are placed, one must place the remaining rooks in the $x\times (n-k)$ subrectangle consisting of those columns of the original rectangle not used for the rooks on $B$.  Placing these rooks from left to right, there will be $x$ choices for the position of the first rook, then $x-m$ choices for the next, and so on, for a total of $x\da_{n-k,m}$ choices.  Thus the sum side is $r_{n,m}(B_x)$ as desired.

On the product side, it will be convenient to consider placing rooks on $B_x$ zone by zone from left to right.  So suppose $z=[s,t]$ is a zone and all rooks in zones to its left have been placed.  Because $z$ is a zone we have 
$\lf b_s \rf_m =\dots = \lf b_t \rf_m = cm$ for some constant $c$.  Also, among all the rooks placed in the columns of $z$, there is at most one which is in the set of squares $\cR$ corresponding to $\rho_m(z)$.  If there are no rooks in $\cR$ then they all go in a rectangle of height $x+cm$.  Thus, using the same ideas as in the previous paragraph, the number of placements is
\beq
\label{0}
(x+cm-(s-1)m)(x+cm-sm)\dots(x+cm-(t-1)m).
\eeq
When there is one rook in $\cR$, say it is in the column with index $i$.  So there are $\rho_m(b_i)$ choices for the placement of this rook and the rest of the rooks go in a rectangle of height $x+cm$.  This gives a count of
\beq
\label{1}
\rho_m(b_i)(x+cm-(s-1)m)(x+cm-sm) \dots(x+cm-(t-2)m).
\eeq
Adding together the contributions from~\ree{0} and~\ree{1} and factoring, we see that the total number of placements is
$$
(x+cm-(s-1)m)\dots(x+cm-(t-2)m)(x+cm-(t-1)m+\rho_m(z)).
$$
Remembering that $\lf b_s \rf_m =\dots = \lf b_t \rf_m = cm$, we see that this is exactly the contribution needed for the product.
\eprf

We should show why our result implies the theorems of Goldman-Joichi-White and Briggs-Remmel.  In both cases, it is clear that the sum sides correspond, so we will concentrate on the products.

For Theorem~\ref{gjw} we take $m=1$.  Since $\lf n \rf_1 = n$ for any $n$, $\rho_1(z) =0$ for any zone $z$ and the two cases in Theorem~\ref{mft} are the same.  So the contribution of the $j$th column to the product is
$$
x+\lf b_j \rf_1 - (j-1)\cdot 1 = x+b_j-j+1
$$
in agreement with the Factorization Theorem.

As far as Theorem~\ref{br},  suppose that $B$ is a singleton board and consider any zone $z=[s,t]$.  If $s\le j<t$ then $b_j$ is a multiple of $m$ and
$$
x+\lf b_j \rf_m - (j-1) m =x +b_j-(j-1)m.
$$ 
And if $j=t$ then $\rho_m(z)=\rho_m(b_t)$  so that
$$
x+\lf b_t\rf_m - (t-1) m +\rho_m(z)   = x + b_t - (t-1)m.
$$
So in either case one gets the same factor as in the Briggs-Remmel result.

%%%%%%%%%%%%%%%%%%%%%%%%%%%%%%%%
%
% 	A $p,q$-ANALOGUE	
%
%%%%%%%%%%%%%%%%%%%%%%%%%%%%%%%%

\section{A $p,q$-analogue}
\label{pqa}

In this section we will derive a $p,q$-analogue of Theorem~\ref{mft}.  Such an analogue was given by Remmel and Wachs~\cite{rw:rtg} for Theorem~\ref{gjw} and was generalized to singleton boards by Briggs and Remmel~\cite{br:mrn}.  Before proving our result, we would like to motivate the statistics we will be using on $p$ and $q$.  

Let $\fS_n$ denote the symmetric group of all permutations $\pi=a_1\dots a_n$ of $\{1,\ldots,n\}$ written in 1-line notation.  One of the most famous statistics on $\fS_n$ is the \emph{inversion number}  which is defined by
$$
\inv \pi = |\{(i,j)\ :\ \text{$i<j$ and $a_i>a_j$}\}|.
$$
By way of illustration, if $\pi=4132$ then $\inv\pi = 4$ because of the pairs $(i,j)=(1,2)$, $(1,3)$, $(1,4)$, and $(3,4)$.
A good source of information about inversions and other statistics is Stanley's book~\cite{sta:ec1}.  Similarly we can consider the \emph{coinversion number} defined by
$$
\coinv \pi =|\{(i,j)\ :\ \text{$i<j$ and $a_i<a_j$}\}|.
$$

To obtain the generating function for these two statistics we define, for any complex number $x$, the \emph{$p,q$-analogue of $x$} to be
$$
[x]_{p,q} = \frac{p^x-q^x}{p-q}.
$$
We will sometimes drop the subscripts if no confusion will result.  Note that if $n$ is a nonnegative integer then, by division, we have
\beq
\label{[n]}
[n]_{p,q} = p^{n-1} + p^{n-2}q + \dots + q^{n-1}.
\eeq
For such $n$, we also define the \emph{$p,q$-factorial}  by
$$
[n]_{p,q}!= [1]_{p,q} [2]_{p,q} \dots [n]_{p,q}.
$$
The following is a famous result of Rodrigues.
\bth[\cite{rod:idp}]
For any integer $n\ge0$,

\vs{5pt}

\eqqed{
\sum_{\pi\in\fS_n} p^{\coinv \pi} q^{\inv \pi} = [n]_{p,q}!
}
\eth

To relate inversions and coinversions to rook placements, it is useful to use the notion of a cohook.  
Let $B$ be a board and $(i,j)\in B$.  The \emph{cohook} of $(i,j)$ is
$$
cH_{i,j} = \{(i',j')\in B\ :\ \text{either $i'=i$ and $j'\le j$ or $j'=j$ and $i'\ge i$}\}.
$$
In Figure~\ref{hd}, the board on the left shows the cohook $cH_{2,3}$ in the Ferrers board $(2,4,4,4)$ which is the set of boxes indicated by the dashed lines.
Hooks and cohooks play an important role in enumerative combinatorics and the representation theory of $\fS_n$.  The reader can consult Sagan's book~\cite{sag:sym} for details.  It will be convenient for us to distinguish two subsets of $cH_{i,j}$, namely the \emph{coarm}
$$
cA_{i,j} = \{(i,j')\in B\ :\  j'<j\}
$$
and the \emph{coleg}
$$
cL_{i,j}= \{(i',j)\in B\ :\  i'> i\}.
$$
So we have the disjoint union $cH_{i,j}=\{(i,j)\}\uplus cA_{i,j} \uplus cL_{i,j}$.

Any permutation $\pi\in\fS_n$ can be considered as a placement $P(\pi)$ of $n$ rooks on the $n\times n$ board $B_n$ where the rooks correspond to the ones in the permutation matrix of $\pi=a_1\dots a_n$.   So the rooks in $P(\pi)$ are in positions
$(i,n-a_i+1)$, $1\le i\le n$, where we complement the second component because we are using Cartesian, rather than matrix, coordinates.  The center board in Figure~\ref{hd} shows the placement for the permutation $\pi=4132$ considered previously.  The \emph{diagram} of $\pi$ is the set of squares 
$$
D(\pi)= B_n\setm (\cup_{(i,j)\in P(\pi)}\  cH_{i,j}).
$$
For $\pi=4132$, $D(\pi)$ is the set of squares in the center board of Figure~\ref{hd} which do not contain dashes.
The diagram of a permutation has applications in Schubert calculus.  See, for example, the article of Bergeron~\cite{ber:ccs}.
It is easy to see that $|D(\pi)|=\inv\pi$.  In fact, we have $i<j$ and $a_i>a_j$ in $\pi$ if and only if $(i,n-a_j+1)\in D(\pi)$.  To obtain $\coinv\pi$ one counts the squares remaining in $B_n$ after removing the $(i,j)$ corresponding to $\pi$, their colegs, and their \emph{arms} which are all squares $(i,j')$ with $j'>j$.

To define an \emph{$m$-cohook} of a square $(i,j)\in B$, one just replaces the coleg of $(i,j)$ with its corresponding level, giving
$$
cH_{i,j;m} = \{(i,j)\}\uplus cA_{i,j} \uplus \{(i',j')\in B\  :\ \text{$i'>i$ and $j,j'$ are in the same level}\}
$$
where the third set of  elements in the disjoint union is called the \emph{$m$-coleg}.
The board on the right in Figure~\ref{hd} displays $cH_{2,3;2}$ for $B=(2,4,4,4)$.  The \emph{$m$-diagram} and \emph{$m$-inversion number} of a placement 
$P\sbe B$ and now defined in the expected way:
$$
D_m(P)=B\setm (\cup_{(i,j)\in P}\ cH_{i,j;m})
\qmq{and}
\inv_m P = |D_m(P)|.
$$
One similarly defines the $m$-coinversion number, $\coinv_m P$.

\bfi
\begin{tikzpicture}
\draw (0,0)--(4,0) (0,1)--(4,1) (0,2)--(4,2) (1,3)--(4,3) (1,4)--(4,4)
(0,0)--(0,2) (1,0)--(1,4) (2,0)--(2,4) (3,0)--(3,4) (4,0)--(4,4)
;
\fill(1.5,2.5) circle(.1);
\draw[very thick, dashed] (1.5,2.5)--(4,2.5) (1.5,2.5)--(1.5,0);
\draw(0,2.5) node {$cH_{2,3}=$};
\end{tikzpicture}
\hs{45pt}
\begin{tikzpicture}
\draw(0,0) grid (4,4);
\draw(.5,.5) node {\scalebox{2}{\symrook}};
\draw(1.5,3.5) node {\scalebox{2}{\symrook}};
\draw(2.5,1.5) node {\scalebox{2}{\symrook}};
\draw(3.5,2.5) node {\scalebox{2}{\symrook}};
\draw[very thick, dashed](.65,.5)--(4,.5) (1.65,3.5)--(4,3.5) (2.65,1.5)--(4,1.5) (3.65,2.5)--(4,2.5)
(.5,.2)--(.5,0) (1.5,3.2)--(1.5,0) (2.5,1.2)--(2.5,0) (3.5,2.2)--(3.5,0)
;
\end{tikzpicture}
\hs{45pt}
\begin{tikzpicture}
\draw (0,0)--(4,0) (0,1)--(4,1) (0,2)--(4,2) (1,3)--(4,3) (1,4)--(4,4)
(0,0)--(0,2) (1,0)--(1,4) (2,0)--(2,4) (3,0)--(3,4) (4,0)--(4,4)
;
\fill(1.5,2.5) circle(.1);
\draw[very thick,dashed] (1.5,2.5)--(4,2.5) (1.5,2.5)--(1.5,0) (2,3.5)--(4,3.5);
\draw[very thick] (0,0)--(4,0) (0,2)--(4,2) (1,4)--(4,4);
\draw(0,2.5) node {$cH_{2,3;2}=$};
\end{tikzpicture}

\capt{Cohooks and the diagram \label{hd}}
\efi

To describe the statistics we will work with, it will be convenient to partition the $m$-diagram as follows.  Given a board $B$ and a placement $P\sbe B$ we let
\begin{align*}
\al_m(P) &= \text{the number of $(i,j)\in B$ above a rook of $P$ and not in the $m$-coleg of any rook in $P$},\\
\be_m(P)&= \text{the number of $(i,j)\in B$ below a rook of $P$ and not in the $m$-coleg of any rook in $P$},\\
\ep_m(P)&= \text{the number of $(i,j)\in B$ in a column with no rook of $P$ and not in the $m$-coleg of any}\\
	  &\hs{15pt}\text{rook in $P$},
\end{align*}
where $\al$, $\be$, and $\ep$ stand for above, below, and empty, respectively.  An example for a $3$-level rook placement on the board $(1,1,2,3,5,7)$ is shown in Figure~\ref{abc} where the boxes for each statistic are labeled with the corresponding Greek letter.
It is clear from the definitions that $\inv_m(P)=\al_m(P)+\ep_m(P)$ and $\coinv_m(P)=\be_m(P)+\ep_m(P)$. 

\bfi
\begin{tikzpicture}
\draw (0,0) -- (6,0) (0,1)--(6,1) (2,2)--(6,2) (3,3)--(6,3) (4,4)--(6,4) (4,5)--(6,5)  (5,7)--(6,7)
(0,0)--(0,1) (1,0)--(1,1) (2,0)--(2,2) (3,0)--(3,3) (4,0)--(4,5) (5,0)--(5,7) (6,0)--(6,7)
;
\draw[very thick] (0,0)--(6,0) (3,3)--(6,3) (5,6)--(6,6);
\draw(3.5,1.5) node {\scalebox{2}{\symrook}};
\draw(4.5,3.5) node {\scalebox{2}{\symrook}};
\draw[very thick, dashed] (3.65,1.5)--(6,1.5) (4.65,3.5)--(6,3.5) (4,.5)--(6,.5) (4,2.5)--(6,2.5) (5,4.5)--(6,4.5) (5,5.5)--(6,5.5)
;
\draw(.5,.5) node {\scalebox{1.5}{$\ep$}};
\draw(1.5,.5) node {\scalebox{1.5}{$\ep$}};
\draw(2.5,.5) node {\scalebox{1.5}{$\ep$}};
\draw(2.5,1.5) node {\scalebox{1.5}{$\ep$}};
\draw(3.5,.5) node {\scalebox{1.5}{$\be$}};
\draw(3.5,2.5) node {\scalebox{1.5}{$\al$}};
\draw(4.5,4.5) node {\scalebox{1.5}{$\al$}};
\draw(5.5,6.5) node {\scalebox{1.5}{$\ep$}};
\end{tikzpicture}
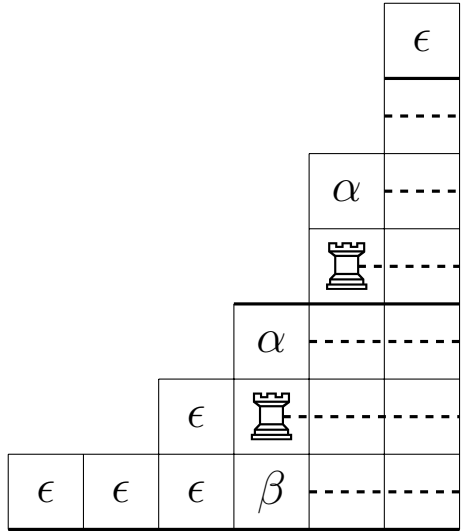
\capt{The $\al$, $\be$ and $\ep$ for a $3$-level placement on $(1,1,2,3,5,7)$ \label{abc}}
\efi

We can now define the $p,q$-analogue of the rook numbers which 
we will be using.  Namely, as in~\cite{br:mrn}, for any board $B$, we let
\beq
\label{r_km[B]}
r_{k,m}[B] = \sum_P p^{\be_m(P)-(c_1+\dots+c_k)m} q^{\al_m(P)+\ep_m(P)}
\eeq
where the sum is over all $m$-level rook placements  $P\sbe B$ with $k$ rooks and $c_1<\dots<c_k$ are the indices of the columns in which the rooks are placed.  The placement in Figure~\ref{abc} would contribute a term of $p^{-26} q^7$ to the sum.  The reason for the unexpected exponent on $p$ is because, when we mimic the proof of Theorem~\ref{mft}, the rooks on the rectangular portion of $B_x$ will never contribute to the $m$-coinversion number of $B$.

We will now generalize the $p,q$-analogue of Theorem~\ref{br} given by Briggs and Remmel to all Ferrers boards.  To state their result, we will use the $p,q$ version of the $m$-falling factorial given by
$$
[x]\da_{n,m}=[x][x-m][x-2m]\dots[x-(n-1)m].
$$
Note that the subscripts no longer refer to $p$ and $q$.
\bth[\cite{br:mrn}]
\label{pqbr}
If  $B=(b_1,\dots,b_n)$  is a singleton board then

\vs{10pt}

\eqqed{
\sum_{k=0}^n p^{xk+m\binom{k+1}{2}}\  r_{k,m}[B]\  [x]\da_{n-k,m} = \prod_{j=1}^n [x+b_j-(j-1)m].
}
\eth

Our generalization will need a finer invariant than the remainder of a zone as used in the $m$-Factorization Theorem.  In particular, if $z=[s,t]$ is a zone of a board $(b_1,\dots,b_n)$ and $i\in z$, then the \emph{$i$th partial remainder} of $z$ is defined by
$$
\rho_{i,m}(z) = \sum_{s\le j\le i} \rho_m(b_j).
$$
We also permit $i=s-1$ in which case $\rho_{i,m}(z)=0$.
\bth
\label{pqmft}
If $B=(b_1,\dots,b_n)$ is any Ferrers board then
$$
\barr{l}
\dil\sum_{k=0}^n p^{xk+m\binom{k+1}{2}}\  r_{k,m}[B]\  [x]\da_{n-k,m}\\[20pt]
= \dil\prod_{j=1}^n 
\begin{cases}
q^{\rho_m(z)}[x+\lf b_j \rf_m -(j-1)m] + \sum\limits_{i\in z} p^{x+\lf b_i \rf_m -(i-1)m}q^{\rho_{i-1,m}(z)}[\rho_m(b_i)] 
&\text{if $j$ is last in zone $z$,}\\[10pt]
[x+\lf b_j \rf_m -(j-1)m] &\text{otherwise.}
\end{cases}
\earr
$$
\eth
\bprf
We will follow the method of proof and notation for Theorem~\ref{mft} while keeping track of the two statistics.  We wish to show that both sides of the equality are the generating function 
\beq
\label{pqB_x}
\sum_{P_x\sbe B_x}  p^{\be_m(P_x)} q^{\al_m(P_x)}
\eeq
where the sum is over all $m$-level placements of $n$ rooks on $B_x$.  Note that since there is a rook in every column of $B_x$, $\al_m(P_x)=\inv_m P_x$ and $\be_m(P_x)=\coinv_m P_x$.

We start with the sum side in the theorem.  Consider all $m$-level placements $P_x\sbe B_x$ with $k$ rooks on $B$ in columns $c_1<\dots<c_k$.  We wish to show that these $P_x$ contribute the same monomials to the theorem's sum and to~\ree{pqB_x}.   We will use the notation $P=P_x\cap B$ and $R$ for the $x\times n$ rectangular portion of $B_x$.

We first claim that  
\beq
\label{pr[B]}
p^{xk+m\binom{k+1}{2}} r_{k,m}[B]
\eeq
gives the contribution to~\ree{pqB_x} of all cells of  $B_x$ which are either in $B$ or are in $R$ below a rook of $B$.  Recalling equation~\ree{r_km[B]}, we see that the monomials in~\ree{pr[B]} are of the form
$$
p^{\be_m(P)+(x-(c_1-1)m)+(x-(c_2-2)m)+\dots+(x-(c_k-k)m)}q^{\al_m(P)+\ep_m(P)}.
$$
For the exponent on $q$, $\al_m(P)$ counts the cells for $\al_m(P_x)$ which are above a rook in $P$ and $\ep_m(P)$ does the same for cells of $B$ which are above a rook in $R$.  For $p$'s exponent, $\be_m(P)$ counts the cells for $\be_m(P_x)$ which are below a rook in $P$ and also in $B$.  The sum in the exponent accounts for the cells for $\be_m(P_x)$ which are below a rook in $P$ and in the rectangle $R$ since there will be $(c_i-i)m$ rows that will be eliminated in column $i$ because of rooks in $R$ to the left of the column.

There remains to show that $[x]\da_{n-k,m}$ is the contribution to~\ree{pqB_x} of all ways to place rooks in the remaining columns of $R$, where we only need consider the cells in $R$ itself.  In the first such column, there are $x$ ways to place a rook and these give exactly the $x$ terms in the sum for $[x]$ in equation~\ree{[n]}.  Similarly, the contribution of the next column is $[x-m]$, and so on, for a total of $[x]\da_{n-k,m}$.  This finishes the proof for the sum side of the theorem.

For the product, as in the proof of the $m$-Factorization Theorem, we break the demonstration up into two cases depending on whether a rook is placed in the region $\cR$ corresponding to the remainder of zone $z$ or not.  If there are no rooks in $\cR$ then, by arguments similar to those already given, the contribution to~\ree{pqB_x} will be
$$
q^{\rho_m(z)}[x+cm-(s-1)m][x+cm-sm]\dots[x+cm-(t-1)m].
$$
When there is a rook in column $i$ of $\cR$, then  $p^{x+\lf b_i \rf_m -(i-1)m}[\rho_m(b_i)] $ gives the contribution of the cells of $B_x$ in column $i$ for all possible placements in that column of $\cR$.  A factor  $q^{\rho_{i-1,m}(z)}$ will come from the fact that all the cells in $\cR$ to the left of column $i$ will be above a rook of $B_x$.  And the cells in $\cR$ to the right of column $i$ contribute nothing since they are in the $m$-coleg of the rook in column $i$.  Considering the contributions from the other columns of $z$ in this case gives a total contribution of 
$$
 p^{x+\lf b_i \rf_m -(i-1)m}q^{\rho_{i-1,m}(z)}[\rho_m(b_i)] [x+cm-(s-1)m][x+cm-sm]\dots[x+cm-(t-2)m].
$$
Summing everything together gives the factors corresponding to $z$ in the product of the theorem.  Multiplying over all zones completes the proof.
\eprf

%%%%%%%%%%%%%%%%%%%%%%%%%%%%%%%%
%
% 	ROOK EQUIVALENCE
%
%%%%%%%%%%%%%%%%%%%%%%%%%%%%%%%%

\section{Rook equivalence}
\label{re}

Two Ferrers boards $B,B'$ are  \emph{$m$-level rook equivalent} if $r_{k,m}(B)=r_{k,m}(B')$ for all $k\ge0$.  In this case we will write $B\equiv_m B'$.  We will drop the ``$m$" and just say ``rook equivalent" if $m=1$.  Foata and Sch\"utzenberger~\cite{fs:rpf} proved a beautiful theorem giving a distinguished board in each equivalence class.  Call a Ferrers board $B=(b_1,\dots,b_n)$ 
\emph{ increasing} if $0<b_1<\dots<b_n$.
\bth[\cite{fs:rpf}]
\label{fs}
Every Ferrers board is rook equivalent to a unique  increasing board.\hqed
\eth

The purpose of this section is to extend Theorem~\ref{fs} to arbitrary $m$.  The Foata-Sch\"utzenberger result was reproved by Goldman-Joichi-White using their Factorization Theorem.  To see the connection, suppose that $B=(b_1,\dots,b_n)$ and 
$B'=(b_1',\dots,b_n')$ are two Ferrers boards.  Although we are writing the boards with the same number of columns, this is no restriction since we can always pad the shorter board with columns of height $0$ on the left.  So $B$ and $B'$ are rook equivalent if and only if they have the same generating function in the falling factorial basis.  By Theorem~\ref{gjw}, this happens if and only if the two vectors
$$
\ze(B)=(-b_1,1-b_2,2-b_3,\dots,(n-1)-b_n) \qmq{and} \ze(B')=(-b_1',1-b_2',2-b_3',\dots,(n-1)-b_n')
$$
are rearrangements of each other since these are the zeros of the corresponding products.  We call $\ze(B)$ the \emph{root vector} of $B$.  For example, if $B=(1,1,3)$ and $B'=(2,3)$ then, rewriting $B'=(0,2,3)$, we have $\ze(B) = (-1,0,-1)$ and $\ze(B')=(0,-1,-1)$ and so $B\equiv B'$.   We should note that  padding a board $B$ with zeros will change the entries of $\ze(B)$.  Also, if $\ze(B)$ is a rearrangement of $\ze(B')$ then the same will be true when padding both $B$ and $B'$ with any given number of zeros.

We now return to considering general $m$.  Define the \emph{$m$-level root vector}  of $B=(b_1,\dots,b_n)$ to be 
$\ze_m(B)=(a_1,\dots,a_n)$ where
$$
a_j =
\begin{cases}
(j-1)m - \lf b_j \rf_m -  \rho_m(z) &\text{if $j$ is the last index in its zone $z$,}\\
(j-1)m - \lf b_j \rf_m  &\text{otherwise.}
\end{cases}
$$
The next result is immediate from Theorem~\ref{mft}.
\bpr
\label{ze}
Ferrers boards $B$ and $B'$ satisfy $B\equiv_m B'$ if and only if $\ze_m(B)$ is a rearrangement of $\ze_m(B')$.\hqed
\epr

\bfi
\begin{tikzpicture}
\foreach \y in {0,1} 
   \draw (0,\y)--(3,\y);
\foreach \y in {2,3,} 
   \draw (1,\y)--(3,\y);
\foreach \x in {1,2,3}
   \draw (\x,0)--(\x,3);
\draw (0,0)--(0,1);
\draw[very thick] (0,0)--(3,0) (1,2)--(3,2);
\draw(-.5,1.5) node {$B=$};
\end{tikzpicture}
\hs{100pt}
\begin{tikzpicture}
\draw(0,0) grid (2,2);
\draw(1,2) grid (2,5);
\draw[very thick] (0,0)--(2,0) (0,2)--(2,2) (1,4)--(2,4);
\draw(-1,1.5) node {$l(B)=$};
\end{tikzpicture}
\capt{A board $B$ and $l(B)$ when $m=2$  \label{l}}
\efi

Our first order of business will be to restrict the class representative problem to considering singleton boards 
$B=(b_1,\dots,b_n)$ since then the $a_j$ are simpler to compute.   Indeed, the argument in the  last paragraph of Section~\ref{sec:mft} shows that in this case
$$
a_j= (j-1)m -b_j
$$
for all $j$.   To describe a singleton board in each equivalence class, let $\cQ_i$ denote the set of squares in the $i$th level of  the quadrant $\cQ$ and, for any board $B=(b_1,\dots,b_n)$, let $l_i=|B\cap \cQ_i|$.  For example, if $m=2$ and $B=(1,3,3)$ as shown on the left in Figure~\ref{l}, then $l_1=5$, $l_2=2$, and $l_i=0$ for $i\ge3$.  For any Ferrers board $B$, if $t$ is the largest index with $l_t\neq 0$ then we let $l(B)=(l_t,l_{t-1},\dots,l_1)$.  We call this function the \emph{$l$-operator} on boards. The board on the right in Figure~\ref{l} shows $l(1,3,3)=(2,5)$.
\ble
\label{l(B)}
For any $m$ and any Ferrers board $B=(b_1,\dots,b_n)$ the sequence $l(B)=(l_t,\dots,l_1)$ is a partition, the Ferrers board $l(B)$ is singleton, and $B\equiv_m l(B)$.
\ele
\bprf
To see that  $l(B)$ is a partition first note that, for any $i>1$, the set of columns of $B\cap \cQ_i$ is a subset of the columns of 
$B\cap \cQ_{i-1}$.  Furthermore, for each of the former columns we have $m$ squares of that column in $B\cap \cQ_{i-1}$.   It follows that $l_i(B)\le l_{i-1}(B)$ as desired.

To show that $l(B)$ is singleton, it suffices to show that if $l_i$ is not a multiple of $m$ then $\lf l_i \rf_m < \lf l_{i-1} \rf_m$.  Let $c$ be the number of columns of $B\cap \cQ_i$ which contain $m$ squares and $d$ be the number of columns containing fewer than $m$ squares.  By assumption $d>0$.  Since these $c+d$ columns of $B\cap \cQ_{i-1}$ must all contain $m$ cells we have
$l_i <(c+d)m\le l_{i-1}$.  Taking floors finishes this part of the proof.

To prove rook equivalence, pick any $I=\{i_1>\dots>i_k\}\sbe\{1,\dots,t\}$.  It suffices to show that the number of ways to place rooks on the levels of $B$ indexed by $I$ is the same as the number of ways to place rooks in the columns of $l(B)$ indexed by $I$.  For the former, if one places the rooks level by level from top to bottom then, since each rook in a higher level rules out  $m$ squares in each level below it, we obtain a count of
$$
l_{i_1}(l_{i_2}-m)(l_{i_3}-2m)\dots(l_{i_k}-(k-1)m).
$$
Now consider  placing the rooks in $l(B)$  column by column from left to right.  Since $l(B)$ is singleton, each rook placed eliminates $m$ squares in each column to its right from consideration.  Thus we obtain the same count as before.
\eprf

The possible $\ze_m$-vectors of singleton boards are easy to characterize.  
\bpr
\label{zequiv}
Consider a vector $\ze=(a_1,\dots,a_n)$ where $a_1=0$.  Let $B=(b_1,\dots,b_n)$ where we define $b_j=(j-1)m-a_j$ for all $j$.  We have that $\ze=\ze_m(B)$  where $B$ is a singleton board if and only if the following conditions are satisfied:
\ben
\item[(i)] $a_{j+1}\le a_j+m$ for $j\ge1$, and
\item[(ii)]  if neither of $a_{j+1}, a_j$ are multiples of $m$ then $\lf a_{j+1} \rf_m \le \lf a_j\rf_m$.
\een
\epr
\bprf
We claim that (i)  and the fact that $a_1=0$  are equivalent to $B$ being a weakly increasing sequence of nonnegative integers.  Since $a_1=0$ we have $b_1=-a_1=0$ which is nonnegative.  And for $j\ge1$
\beq
\label{diff}
b_{j+1}-b_j=(jm-a_{j+1})-((j-1)m-a_j)=m+a_j-a_{j+1}.
\eeq
So (i) is equivalent to $B$ being weakly increasing.  A similar argument shows that (ii) is equivalent to the singleton condition.
\eprf

We are finally in a position to give distinguished representatives for the $m$-level equivalence classes. A  Ferrers board 
$B=(b_1,\dots,b_n)$ is called \emph{$m$-increasing} if $b_1>0$ and $b_{j+1}\ge b_j+m$ for $j\ge1$.  Note that a $1$-increasing board is increasing in the sense of the definition before Theorem~\ref{fs}.  Also note that, although  most properties of Ferrers boards are not affected by padding with columns of length zero, the $m$-increasing condition will be destroyed.
\bth
\label{minc}
Every Ferrers board is $m$-level rook equivalent to a unique $m$-increasing board.
\eth
\bprf
Clearly any $m$-increasing board is a singleton board.  So, by Lemma~\ref{l(B)}, it suffices to show that any singleton board $B$ is $m$-rook equivalent to a unique $m$-increasing board.  An example of the construction of this board is given after this proof to illustrate the technique.  Let $N=|B|+1$ and pad $B$ with columns of zeros so as to write $B=(b_1,\dots,b_N)$.  By the choice of $N$, any board $B'$ which is $m$-equivalent to $B$ can be written as 
$B'=(b_1',\dots,b_N')$  and $b_1=b_1'=0$.  Letting $\ze_m(B)=(a_1,\dots,a_N)$ and $\ze_m(B')=(a_1',\dots,a_N')$, the choice of $N$ also ensures that  $a_1=a_1'=0$ and $a_i, a_i'\ge 0$ for all $i$.  

We claim  that a singleton board $B'=(b_1',\dots,b_n')$ will be $m$-increasing if and only $\ze_m(B')=(a_1',\dots,a_n')$ is weakly decreasing.  Indeed, this follows directly from~\ree{diff}.

We first show existence.  By the previous paragraph, we wish to rearrange $\ze_m(B)$ in such a way that the portion of the rearrangement corresponding to nonzero entries of the board is weakly decreasing.  And the portion of the rearrangement corresponding to zero entries of the board must be of the form $0,m,2m,\dots,cm$ for some $c$.  So choose $cm$ to be the largest multiple of $m$ in $\ze=\ze_m(B)$.  We claim that the elements $0,m,\dots,(c-1)m$ also occur in $\ze$.
Consider the multiples of $m$ in the first zone of $B$.  Since $b_1=0$, these will all be zeros.  It follows that in $\ze$ we have entries $0,m,\dots,im$ for some $i$.  Now consider the multiples of $m$ in the next zone of $B$.  Since they are all at least $m$ and there is at most one non-multiple of $m$ in the first zone, in $\ze$ we will have entries $jm,(j+1)m,\dots$ for some $j\le i+1$.  This will continue until we get to a zone giving rise to the entry $cm$. The inequalities on the factors of $m$ between adjacent zones  imply the claim.

We now define $\ze'=(a_1',\dots,a_N')$ 
where 
$$
(a_1',a_2',a_3',\dots,a_{c+1}')=(0,m,2m,\dots,cm)
$$
and $(a_{c+2}',a_{c+3}',\dots,a_N')$ is the rest of $\ze$ arranged in weakly decreasing order.  Since $a_1'=0$, we can show that $\ze'$ corresponds to an $m$-increasing board by checking conditions (i) and (ii) of Proposition~\ref{zequiv}.  Condition (i) is clearly true for 
$j\le c$.  For $j=c$, one can show, by using a proof as in the previous paragraph and the choice of $cm$ as the largest multiple of $m$ in $\ze$, that $a_{c+2}'<(c+1)m$.  So (i) also holds in this case.  And for $j>c$ the fact that the sequence is weakly decreasing makes (i) a triviality.  The same reasoning as for (i) shows that (ii) must also hold. Thus, defining $B'=(b_1',\dots,b_N')$ where $b_j'=(j-1)m-a_j'$ for all $j$ results in a singleton board.  Furthermore, by construction, $b_1'=\dots=b_{c+1}'=0$ and 
$(b_{c+2}',\dots,b_N')$ is $m$-increasing.  Hence removing the zeros from $B'$ leaves the desired $m$-increasing board.

To show uniqueness, suppose $\ze'=(a_1',\dots,a_N')$ is a rearrangement of $\ze$ corresponding to a padded $m$-increasing board.  Then $\ze'$ must start with $0,m,\dots,cm$ for some $c$ and be weakly decreasing thereafter by equation~\ree{diff}.  Without loss of generality, we can assume $a_{c+2}'\neq (c+1)m$, since if the two are equal we can just add $a_{c+2}'$ to the initial run of multiples of $m$.  So $\ze'$ will be the rearrangement of $\ze$ constructed in the existence proof as long as $cm$ is the largest multiple of $m$ in $\ze$.  But if $cm$ is not the largest multiple of $m$ in $\ze$ then $(a_{c+2}',\dots,a_N')$ must contain an element greater than or equal to $(c+1)m$.  And since this portion of $\ze'$ is weakly decreasing, this forces $a_{c+2}'>(c+1)m$ since equality was ruled out earlier.  But then $a_{c+1}'=cm$ and $a_{c+2}'$ do  not satisfy condition (i) of Proposition~\ref{zequiv}, contradicting the fact that $\ze'$ corresponds to a singleton board.  This finishes the proof of uniqueness and of the theorem.
\eprf

To illustrate the construction in the previous proof, consider $m=2$ and the singleton board $B=(1,2,2,3)$.  Now $N=1+2+2+3=8$ and we pad $B$ with zeros to length $8+1=9$, obtaining $B=(0,0,0,0,0,1,2,2,3)$.  Thus $\ze=\ze_2(B)=(0,2,4,6,8,9,10,12,13)$.  The largest multiple of $2$ in $\ze$ is $12$, so we rearrange $\ze$ to begin with the multiples of $2$ up through $12$ and then decrease.  The result is $\ze'=(0,2,4,6,8,10,12,13,9)$ with associated board $B'=(0,0,0,0,0,0,0,1,7)$.  Removing the initial zeros, we get the $2$-increasing board $(1,7)$ which is easily seen to be $2$-level rook equivalent to $B$.

%%%%%%%%%%%%%%%%%%%%%%%%%%%%%%%%
%
% 	ENUMERATION OF SINGLETON BOARDS
%
%%%%%%%%%%%%%%%%%%%%%%%%%%%%%%%%

\section{Enumeration of singleton boards}
\label{esb}

Goldman, Joichi, and White~\cite{gjw:rtI} used their factorization theorem to give a simple formula for the size of a given rook equivalence class.  The basic idea is to count, for any board $B$,  the number of rearrangements of $\ze_1(B)$ which correspond to a Ferrers board.  To state their result, given any finite vector $\nu$ of nonnegative integers, we let $n(\nu)=(n_0,n_1,\dots)$ be defined by
$$
n_i=\text{the number of copies of $i$ in $\nu$.}
$$ 
So $n_i(\nu)=0$ if $i<0$ or $i$ is sufficiently large.

\bth[\cite{gjw:rtI}]
\label{gjwcard}
Let $B=(b_1,\dots,b_N)$ be a Ferrers board where $N=|B|+1$, and suppose $n(\ze_1(B))=(n_0,n_1,\dots)$.  The number of Ferrers boards in the equivalence class of $B$ is

\vs{10pt}

\eqqed{
\dil\prod_{i\ge 1} \binom{n_i+n_{i-1}-1}{n_i}.
}
\eth

Because the entries of $\ze_m(B)$ are more complicated for $m\ge2$, we will not be able to count all boards in an $m$-level equivalence class.  But we can at least enumerate the singleton boards.  The formula will be in terms of multinomial coefficients.
\bth
\label{card}
Let $B=(b_1,\dots,b_N)$ be a singleton board where $N=|B|+1$, and suppose $n(\ze_m(B))=(n_0,n_1,\dots)$.  Then the number of singleton boards which are $m$-level rook equivalent to $B$ is
$$
\prod_{i\ge0} \binom{n_{im}+n_{im+1}+\dots+n_{im+m}-1}{n_{im}-1,n_{im+1},n_{im+2},\dots,n_{im+m}}.
$$
\eth
\bprf
By the choice of $N$ we have that $\ze=\ze_m(B)$ begins with a zero and has all entries nonnegative.  So it suffices to compute the number of rearrangements of $\ze$ beginning with $0$ and satisfying conditions (i) and (ii) of Proposition~\ref{zequiv}.  Let $d$ be the maximum entry of $\ze$.  If $d=0$ then the result is easy to verify, so assume $d>0$.  Let $cm$ be the largest nonnegative multiple of $m$ with $cm<d$.  Note that $cm$ exists since $d>0$ and also that, by the argument given in the proof of Theorem~\ref{minc}, $n_{cm}>0$.

Consider the vector $\ze'=(a_1',\dots,a_n') $ obtained from $\ze$ by removing all entries which are larger than $cm$.  We claim that $\ze'=\ze_m(B')$ for some singleton board $B'$.   As usual, we use Propostion~\ref{zequiv}.  Certainly $a_1'=a_1=0$ since none of the zeros were removed from $\ze$.  Suppose, towards a contradiction, that condition (i) is false in that $a_{j+1}'>a_j'+m$ for some $j$.  Let $a_i$ be the element of $\ze$ corresponding to $a_j'$.  But since we removed the largest elements of $\ze$ we have 
$$
a_{i+1}\ge a_{j+1}' >a_j'+m = a_i+m
$$
which is impossible.  A similar contradiction demonstrates (ii), and our claim is proved.

Now, by induction, it suffices to show that the number of rearrangements of $\ze$ which come from a given $\ze'$ by inserting elements larger than $cm$ is
$$
\binom{n_{cm}+n_{cm+1}+\dots+n_{cm+m}-1}{n_{cm}-1,n_{cm+1},n_{cm+2},\dots,n_{cm+m}}.
$$
Consider  elements $a_j,a_{j+1}$ in $\ze$.  First note that if  $a_j$ comes from $\ze'$ and $a_{j+1}>cm$ then we must have $a_j=cm$.  If this were not the case then, since $a_j<cm<a_{j+1}$, to make condition (i) true neither $a_j$ nor $a_{j+1}$  would be multiplies of $m$.  But by the same pair of inequalities  we would have $\lf a_{j+1} \rf_m \ge cm > \lf a_{j} \rf_m$ which contradicts condition (ii).  Thus we can insert elements greater than $cm$ only after copies of $cm$ itself.   In addition, any $a_{j+1}$ with $cm<a_{j+1}\le cm+m$ can come after $a_j=cm$ as it is easily verified that we always have conditions (i) and, vacuously, (ii) for such $a_{j+1}$.

We also claim that the elements larger than $cm$ can be arranged in any order with respect to each other.  To see this, suppose $cm<a_j,a_{j+1}\le cm+m$.  Condition (i) is immediate because of the given bounds.  And if neither is a multiple of $m$ then we have $cm<a_j,a_{j+1}<cm+m$ which implies condition (ii).

Finally, if $a_j>cm$ and $a_{j+1}$ comes from $\ze'$ then $a_{j+1}<a_j$ and conditions (i) and (ii) are trivial.  So an element greater than $cm$ can be followed by any element of $\ze'$.

We now enumerate the number of $\ze$ coming from $\ze'$ using the structural properties from the previous three paragraphs.
There are $n_{cm}$ copies of $cm$ and $n_{cm+1}+n_{cm+1}+\dots+n_{cm+m}$ elements to be inserted after these copies where the space after a copy can be used multiple times.  And any element of $\ze'$ can follow the inserted elements.  So the total number of choices for this step is the binomial coefficient
$$
\binom{n_{cm}+n_{cm+1}+\dots+n_{cm+m}-1}{n_{cm+1}+n_{cm+2}+\dots+n_{cm+m}}.
$$
Now we need to arrange the elements greater than $cm$ among themselves.  This can be done arbitrarily, so the number of choices is
$$
\binom{n_{cm+1}+n_{cm+2}+\dots+n_{cm+m}}{n_{cm+1},n_{cm+2},\dots,n_{cm+m}}.
$$
Multiplying the two displayed expressions and canceling $(n_{cm+1}+n_{cm+2}+\dots+n_{cm+m})!$ gives the desired result.
\eprf

Note that the result just proved does indeed generalize Theorem~\ref{gjwcard}.  This is because when $m=1$, every board is a singleton board.  And the products clearly conicide in this case.

%%%%%%%%%%%%%%%%%%%%%%%%%%%%%%%%
%
% 	FILE PLACEMENTS	
%
%%%%%%%%%%%%%%%%%%%%%%%%%%%%%%%%

\bfi
\begin{tikzpicture}
\draw(0,0) grid (6,2);
\draw(2,2) grid (6,3);
\draw(.5,.5) node {\scalebox{2}{\symrook}};
\draw(1.5,.5) node {\scalebox{2}{\symrook}};
\draw(3.5,.5) node {\scalebox{2}{\symrook}};
\draw(2.5,2.5) node {\scalebox{2}{\symrook}};
\draw(5.5,2.5) node {\scalebox{2}{\symrook}};
\end{tikzpicture}
\capt{A file placement \label{F}}
\efi

\section{File placements}
\label{fp}

Thus far our focus has been to keep the sum side of  Theorem~\ref{br} the same and modify the product side to get equality for all Ferrers boards.  Another possibility would be to keep the product side the same, expand it in the $m$-falling factorial basis, and see if the coefficients of the linear combination count anything.  This will be our approach in the current section.  The case $m=2$ was considered earlier in a paper of Haglund and Remmel~\cite{hr:rtp}.

It turns out that these coefficients count weighted file placements.  A \emph{file placement} on a board $B$ is $F\sbe B$ such that no two rooks (elements) of $F$ are in the same column.  However, we permit rooks to be in the same row.  
Figure~\ref{F} displays such a placement on the Ferrers board $(2,2,3,3,3,3)$.

We let
$$
f_k(B) =\text{the number of file placements $F\sbe B$ with $k$ rooks.}
$$
It is easy to count such placements.  If $B$ has $b_j$ squares in column $j$ for $1\le j\le n$ ($B$ need not be a Ferrers board) then $f_k(B)=e_k(b_1,\dots,b_n)$ where $e_k$ is the $k$th elementary symmetric function.  So in order to get more interesting results, we will weight file placements.

Fix, as usual, a positive integer $m$.  Given a board $B$ and a file placement $F\sbe B$, let $t$ be the largest index of a row in $B$ and consider $y_1,\dots,y_t$ where $y_i$ is the number of rooks of $F$ in row $i$.  Define the \emph{$m$-weight} of $F$ to be
$$
\wt_m F = 1\da_{y_1,m} 1\da_{y_2,m}\dots 1\da_{y_t,m}.
$$
In the example of Figure~\ref{F} with $m=3$ we have
$$
\wt_m F = 1\da_{3,3} \cdot 1\da_{0,3} \cdot 1\da_{2,3} = (1)(-2)(-5)\cdot (1) \cdot (1)(-2) = -20.
$$
Note that if $F$ is actually a rook placement then $\wt_m F=1$ because $1\da_{0,m}=1\da_{1,m}=1$ for any $m$.

Given a board, $B$, we define the associated \emph{$m$-weighted file placement numbers} to be
$$
f_{k,m}(B) =\sum_F \wt_m F
$$
where the sum is over all file placements $F\sbe B$ with $k$ rooks.  These are the coefficients which we seek.
\bth[$m$-weight Factorization Theorem]
\label{mwft}
For any   Ferrers board $B=(b_1,\dots,b_n)$ 
\beq
\label{f_km}
\sum_{k=0}^n f_{k,m}(B) x\da_{n-k,m} = \prod_{j=1}^n (x+b_j-(j-1)m).
\eeq
\eth
\bprf
In the manner to which we have become accustomed, we consider the board $B_x$ obtained by attaching an $x\times n$ rectangle $R$ to $B$ where $x$ is a nonnegative multiple of $m$.   Consider \emph{mixed placements} $F\sbe B_x$ which are file placements when restricted to $B$, but satisfy the $m$-level condition when restricted to $R$.  We will compute $S=\sum_F \wt_m F$ where the sum is over the mixed placements $F$ on $B_x$ with $n$ rooks.

We first recover the sum side of equation~\ree{f_km}.  The  mixed placements with $k$ rooks on $B$ will contribute  
$f_{k,m}(B)$ to $S$ from these rooks.   And $x\da_{n-k,m}$ will be the contribution from the $n-k$ rooks on $R$ by the remark about the weight for rook placements above.   So
\beq
\label{k}
f_{k,m}(B) x\da_{n-k,m}=\sum_{F^k} \wt_m F^k
\eeq
where the sum is over all mixed placements $F^k\sbe B_x$ having $k$ rooks in $B$.  Now summing over $k$ gives us the desired equality.

To obtain the product, let $B'$  and $B_x'$ be  $B$ and $B_x$ with their $n$th columns removed, respectively.  By induction on $n$, it suffices to prove that
\beq
\label{B'}
\sum_{k=0}^n f_{k,m}(B) x\da_{n-k,m}=(x+b_n-(n-1)m)\sum_{k=0}^{n-1} f_{k,m}(B') x\da_{n-k-1,m}.
\eeq
Comparing equations~\ree{k} and~\ree{B'}, we see it suffices to show that, for any mixed placement $F'\sbe B_x'$,
\beq
\label{F'}
(x+b_n-(n-1)m) \wt_m F' = \sum_F \wt_m F
\eeq
where the sum is over all mixed placements $F\sbe B_x$ whose restriction to $B_x'$ is $F'$.  To this end, let $y_0$ be the number of rooks in $F'$ which are in $R$.  Also let $y_i$, $1\le i\le b_n$, be the number of rooks in the $i$th row of $F'\cap B'$.  Since every column of $F'$ has a rook, we have
\beq
\label{y}
y_0+y_1+\dots+ y_{b_n}=n-1.
\eeq

We now consider two cases depending on whether the rook in column $n$ of $F$ lies in $B$ or in $R$.  If it lies in $R$ then, by the $m$-level condition, there are $x-y_0m$ places for the rook.  Since these placements are in rows not occupied by rooks of $F'$, each of them contributes a factor of $1$ to the weight for a total change in weight of $x-y_0m$ from this case.  Now suppose that the rook lies in $B$, say in the $i$th row.  Then in passing from $F'$ to $F$, the weight is changed from $1\da_{y_i,m}$ to $1\da_{y_i+1,m}$.  This means that the weight is multiplied by $1-y_i m$ when placing a rook in row $i$ of $B$.  Adding together all the contributions and using equation~\ree{y} gives
$$
x-y_0m+\sum_{i=1}^{b_n} (1-y_im)=x+b_n-(n-1)m.
$$
This completes the proof of equation~\ree{F'} and of the theorem.
\eprf

We note that Theorem~\ref{mwft} is another generalization of the Factorization Theorem.  Indeed, when $m=1$, then any file placement having a row with $y\ge2$  rooks will have a factor of $1\da_{y,1}=0$. And any rook placement will have a weight of $1$.  Thus $f_{k,1}(B)=r_{k,1}(B)$.

%%%%%%%%%%%%%%%%%%%%%%%%%%%%%%%%
%
% 	WEIGHT EQUIVALENCE CLASSES 	
%
%%%%%%%%%%%%%%%%%%%%%%%%%%%%%%%%

\section{Weight equivalence classes}
\label{wec}

Given $m$, define  two boards $B,B'$ to be \emph{$m$-weight file equivalent}, written $B\approx_m B'$, if 
$f_{k,m}(B)=f_{k,m}(B')$ for all $k\ge0$. 
Our goal in this section is to find distinguished representatives for the $m$-weight file equivalence classes.  Interestingly, our result will be  dual to the one for $m$-level rook equivalence in the sense that the inequalities will be reversed.  In order to define the representatives, we will have to assume that all our boards start with at least one zero.  So for this section we will write our  Ferrers boards in the form $B=(b_0,b_1,\dots,b_n)$  where $b_0=0$.

   We can use Theorem~\ref{mwft} to test $m$-weight file equivalence.  The \emph{$m$-weight root vector} of a Ferrers board $B=(b_0,b_1,\dots,b_n)$ is
$$
\om_m(B)=(-b_0,m-b_1,2m-b_2,\dots,nm-b_n).
$$
From the $m$-weight Factorization Theorem we immediately get the following.
\bpr
\label{om}
Ferrers boards $B$ and $B'$ satisfy $B\approx_m B'$ if and only if $\om_m(B)$ is a rearrangement of $\om_m(B')$.\hqed
\epr

We will also need a characterization of the vectors which can be $m$-weight root vectors.  The proof of the next result is similar to that of Proposition~\ref{zequiv} and so is omitted.
\bpr
\label{oequiv}
Consider a vector $\om=(a_0,a_1,\dots,a_n)$.   Let $B=(b_0,b_1,\dots,b_n)$ where we define $b_j=jm-a_j$ for all $j$.  
We have that $\om=\om_m(B)$  where $B$ is a Ferrers board  if and only if the following conditions are satisfied:
\ben
\item[(i)] $a_0=0$, and
\item[(ii)] $a_{j+1}\le a_j+m$ for $j\ge0$.\hqed
\een
\epr

Now define a Ferrers board $B=(b_0,b_1,\ldots,b_n)$ to be \emph{$m$-restricted} if $b_{j+1}\le b_j+m$ for all $j\ge0$.  
We now show that the $m$-weight file equivalence class of any Ferrers board contains a unique $m$-restricted board.  An example of the construction of this board follows the proof.
\bth
\label{mres}
Every Ferrers board $B=(b_0,\dots,b_n)$ is $m$-weight file equivalent to a unique $m$-restricted board.
\eth
\bprf
Similarly to the proof  of Theorem~\ref{minc}, we rewrite $B=(b_0,\dots,b_N)$ where $N=|B|$.  This assures us that any equivalent board $B'=(b_0',\dots,b_N')$  has $\om_m(B')$ which is nonnegative and starts with zero.  Also, using equation~\ree{diff}, we see that $B'$ is $m$-restricted if and only if $\om(B')$ is weakly increasing.  So consider 
$\om'=(a_0',\dots,a_N')$ which is the unique weakly increasing rearrangement of $\om=\om_m(B)=(a_0,\dots,a_N)$.    It suffices to show that $\om'=\om_m(B')$ for some board $B'$.  So we just need to check the two conditions of Proposition~\ref{oequiv}.  Condition (i) follows from the choice of $N$ and the fact that $\om$ is nonnegative.  For condition (ii), assume, towards a contradiction, that there is an index $i$ such that $a_{i+1}'>a_i'+m$.  Let $a_j$ be the element of $\om$ which was rearranged to become $a_{i+1}'$.  Then $a_j=a_{i+1}'>0$ and so there must be an element $a_k$ with $k<j$ and $a_k<a_j$.  Let $k$ be the largest such index.  By the choice of $k$ and the fact that $\om$ satisfies (ii), we must have $a_k\ge a_j-m$.  
Thus 
$$
a_i'<a_{i+1}'-m=a_j-m\le a_k< a_j=a_{i+1}'.
$$
But then when rearranging $\om$ in weakly increasing order, $a_k$ should have been placed between $a_i'$ and $a_{i+1}'$, a contradiction. 
\eprf

By way of illustration, suppose that we take $m=2$ and $B=(1,5)$.  This gives $N = |B|=6$ 
and 
$\om=\om_2(B) =(0,2,4,6,8,10,12)-(0,0,0,0,0,1,5) = (0,2,4,6,8,9,7)$.  The weakly increasing rearrangement of $\om$ is
$\om'=(0,2,4,6,7,8,9)$ and so $B'=(0,2,4,6,8,10,12)-(0,2,4,6,7,8,9)=(0,0,0,0,1,2,3)$.

There is a close relationship between the $m$-increasing boards introduced in Section~\ref{re} and $m$-restricted boards.  This is easy to see if $m=1$.  In this case, board $B$ is $1$-increasing if and only if its transpose $B^t$ (obtained by interchanging rows and columns) is $1$-restricted.  Indeed, a Ferrers board is $1$-increasing if and only if 
the northwestern boundary of $B$ contains no horizontal line segment of length at least $2$.  And a board is $1$-restricted if and only if this boundary contains no vertical line segment of length at least $2$.  Note also that when $m=1$, the $l$-operator of Section~\ref{re} satisfies $l(B)=B^t$.  In generalizing these ideas to all $m$, it is the $l$-operator which is key as the next result shows.
\bpr
\label{lprop}
The $l$-operator has the following properties.
\ben
\item[(i)]  If $B$ is $m$-restricted then $l(B)$ is $m$-increasing.
\item[(ii)]  If $B$ is $m$-increasing then $l(B)$ is $m$-restricted.
\item[(iii)]  If $B$ is a singleton board then $l^2(B)=B$, so $l$ is an involution on singleton boards.
\een
\epr
\bprf
(i)  Let $B=(b_1,\dots,b_n)$ be $m$-restricted and $l(B)=(l_t,\ldots,l_1)$.  Keeping in mind that the subscripts in $l(B)$ are decreasing, we wish to show that $l_i\ge l_{i+1}+m$.  Let $B_j$ be the set of cells in column $j$ and let $c_j=|B_j\cap \cQ_i|$ and $d_j=|B_j\cap \cQ_{i+1}|$ for all $j$.  Now $l_i-l_{i+1}=\sum_j (c_j-d_j)$ and $c_j-d_j\ge0$ for all $j$.  Since $B$ is $m$-restricted, there is an index $k$ such that $B_k$ has its highest cell in $\cQ_i$.  Let $k$ be the largest such index.  Using the fact that $B$ is $m$-restricted again forces $B_{k+1}$ to have its highest cell in $\cQ_{i+1}$.  Thus, using the $m$-restricted condition a third time,
$$
l_i-l_{i+1}\ge (c_k-d_k)+(c_{k+1}-d_{k+1})= (b_k-(i-1)m)-0+m-(b_{k+1}-im)=b_k-b_{k+1}+2m\ge m.
$$
which is what we wished to prove

(ii)  This is similar to (i), finding an upper bound for $l_i-l_{i+1}$ using the fact that, for an $m$-increasing board, there is at most one $B_k$ with its highest square in $\cQ_i$ for each $i$.  Details are left to the reader.

(iii)  Induct on $|B|$.  Given $B$, let $B'$ be the board obtained by removing its first level.  Since $B$ is singleton so is $B'$ and thus, by induction, $l^2(B')=B'$.  If $l(B)=(l_t,\dots,l_1)$ then the definition of the $l$-operator shows that $l(B')=(l_t,\dots,l_2)$ and $l_1=|B\cap \cQ_1|$.   Applying $l$ to $l(B)$ we see that the column for $l_1$ in $l(B)$ adds $m$ to every column of $l^2(B')$.  Hence every column of $l^2(B)$ which contains cells in $\cQ_i$ for $i\ge2$ agrees with the corresponding column in $B$.   Also, again by definition of $l$, those columns of $l^2(B)$ which lie wholly in $\cQ_1$ are obtained from the column for $l_1$ in $l(B)$ by breaking it into columns of length $m$ and a column of length $\rho_m(l_1)$.  Since this is also the unique way to complete $l^2(B)$ so that it is a singleton board, it must be that $l^2(B)=B$ as desired.
\eprf

We now return to considering $m$-level rook equivalence classes as in Section~\ref{re}.  Using the proposition just proved, we obtain a second distinguished representative in each $m$-level rook equivalence class.

\bco
Every Ferrers board is $m$-level rook equivalent to a unique $m$-restricted singleton board.
\eco
\bprf
By Theorem~\ref{minc}, we know that each class has a representative $B$ which is $m$-increasing and so also a singleton board.  Applying the previous proposition and Lemma~\ref{l(B)}, we see that $l(B)$ is an $m$-restricted singleton board in the class.  If there is a second such board $B'\neq l(B)$ then,  by Proposition~\ref{lprop} again, $l(B')$ and $l^2(B)=B$ will be distinct $m$-increasing singleton boards in the class, contradicting the uniqueness part of Theorem~\ref{minc}.
\eprf

%%%%%%%%%%%%%%%%%%%%%%%%%%%%%%%%
%
% 	WEIGHT EQUIVALENCE CLASS SIZES	
%
%%%%%%%%%%%%%%%%%%%%%%%%%%%%%%%%

\section{Weight equivalence class sizes}
\label{wecs}

\subsection{Enumeration}
\label{e}

In this subsection we will generalize Theorem~\ref{gjwcard} to $m$-weight equivalence classes.  We will give two proofs of our result:
one from first principles using ideas similar to those in the proof of Theorem~\ref{card} and one using a connection with the theory of $q,t$-Catalan numbers.
\bth
\label{card2}
Let $B=(b_0,\dots,b_N)$ be a Ferrers board where $N=|B|$, and suppose $n(\om_m(B))=(n_0,n_1,\dots)$.  The number of Ferrers boards in the $m$-weight equivalence class of $B$ is
$$
\prod_{i\ge 1} \binom{n_i+n_{i-1}+\dots+n_{i-m}-1}{n_i}.
$$
\eth
\bprf
We use Proposition~\ref{om} and count the number of rearrangements of $\om=\om(B)$ which correspond to a Ferrers board.
Let $d$ be the maximum value of an entry of $\om$.  Our assumptions imply that  $d$ is nonnegative and all entries of $\om$ are between between  $0$ and $d$ inclusive.  Consider $\om'$ which is obtained from $\om_m(B)$ by removing all values equal to $d$.  Using Proposition~\ref{oequiv}, it is easy to see that $\om'=\om_m(B')$ for some Ferrers board $B'$.  By induction, it suffices to show that the number of rearrangements of $\om$ which come from a given $\om'$ is
\beq
\label{binom}
\binom{n_d+n_{d-1}+\dots+n_{d-m}-1}{n_d}.
\eeq
By condition (ii) in the proposition just cited, we can insert $d$ after any element of $\om'$ which is at least $d-m$.  So the number of places for insertion is $n_{d-m}+\dots+n_{d-1}$.  Since we need to insert $n_d$ copies of $d$ and more than one copy can go in a given place, the number of choices  is given by~\ree{binom}.
\eprf

For our second proof of Theorem~\ref{card2}, we will need some background on $q,t$-Catalan numbers.  These polynomials are important in combinatorics and in the study of diagonal harmonics.  See Haglund's book~\cite{hag:qtc} for more information in this regard.  Catalan numbers are intimately connected with lattice paths.  A \emph{$NE$ lattice path} $p$ is a walk on the integer lattice $\bbZ^2$ starting at $(0,0)$ and using unit steps north (parallel to the $y$-axis) and east (parallel to the $x$-axis).  We write $p=s_1,s_2,\dots,s_n$ where each $s_i=N$ or $E$.  The (ordinary) \emph{Catalan numbers} can be defined by
$$
C_n=\text{the number of lattice paths ending at $(n,n)$ which stay weakly below the line $y=x$.}
$$
We will be concerned with $q,t$-analogues of the \emph{$m$-Catalan numbers}, also called \emph{higher Catalan numbers},
$$
C_{n,m}=\text{the number of lattice paths ending at $(n,nm)$ which stay weakly below the line $y=mx$.}
$$

To describe the statistics we will be using as powers of $q$ and $t$, we will express these concepts in terms of Ferrers boards.  We warn the reader that the conventions of people working with the $C_{n,m}(q,t)$ differ from the ones in this paper.  But their diagrams can be obtained from ours by reflecting in the line $y=x$.  Consider the \emph{$m$-triangular board}
$$
\De_{n,m}=(0,m,2m,\dots,(n-1)m).
$$
To make the connection with lattice paths, we make the convention that the right-hand border of $\De_{n,m}$, as well as any other Ferrers board $B$ with $n$ columns, extends up to the point $(n,nm)$.  We call the $NE$ lattice path from $(0,0)$ to $(n,nm)$ consisting of the northwest-most points in $B$ the \emph{boundary} of $B$.  So any $NE$ lattice path to $(n,mn)$ and staying below $y=mx$ must lie weakly southeast of the boundary of $\De_{n,m}$.  The diagram on the left in Figure~\ref{Delta} shows $\De_{4,2}$ with its boundary path $p= E,N,N,E,N,N,E,N,N,E,N,N$ thickened.

Now if $B=(b_0,\ldots,b_{n-1})$ we write $B\sbe\De_{n,m}$ and say that $B$ is \emph{contained in} $\De_{n,m}$ whenever 
$b_j\le jm$ for $0\le j<n$.  We will be interested in
$$
\cB_{n,m}=\text{the set of Ferrers boards $B$ such that $B\sbe\De_{n,m}$}.
$$
In the central diagram of Figure~\ref{Delta} the white squares show $B=(0,0,3,4)$ inside $\De_{4,2}$ and the boundary of $B$ is thickened.  Note that there is a bijection between $\cB_{n,m}$ and the lattice paths counted by $C_{n,m}$ which is given by taking the  boundary of the board.  Note also that if $B=(b_0,\dots,b_{n-1})\sbe\De_{n,m}$ then the $j$th coordinate of $\om_m(B)$,  $jm-b_j$, is exactly the number of squares in the $j$th column of the set difference $\De_{n,m}\setm B$.  Here, and for the rest of this section, we will always start indexing columns at $j=0$ and $b_0=0$ for any board $B$.  In Figure~\ref{Delta}, $\om_2(B)=(0,2,1,2)$ and the squares counted by this vector are shaded.

Our first definition of $C_{n,m}(q,t)$ will be in terms of the $m$-area and $m$-diagonal inversion statistics on boards. 
If $B\sbs\cQ$ is any board then we define its \emph{$m$-area} to be 
$$
\area_m(B)=|B|.
$$
In particular, if
$B\sbe\De_{n,m}$ and $\om_m(B)=(z_0,\dots,z_{n-1})$  then
$$
\area_m  (\De_{n,m}\setm B) =z_0+\dots+z_{n-1}.
$$
To motivate the second statistic, if $\pi=a_1\dots a_n$ is any sequence of real numbers then we can define its \emph{inversion number} in the same way as was done for permutations, namely
$$
\inv\pi=|\{i<j\ :\ a_i-a_j>0\}|.
$$
The \emph{$m$-diagonal inversion number} of $\pi$ is defined by
\beq
\label{dinvdef}
\dinv_m\pi =\sum_{k=0}^{m-1} |\{i<j\ :\ 0\le a_i-a_j+k \le m\}|.
\eeq
Note that, unlike other quantities we have defined, the $1$-diagonal inversion number does not coincide with the inversion number, although they are clearly related.  We extend this definition to boards by setting
$$
\dinv_m B =\dinv_m(\om_m(B)).
$$
With these statistics, we can define the \emph{higher $q,t$-Catalan numbers} by
\beq
\label{cO}
C_{n,m}(q,t)=\sum_{B\in\cB_{n,m}} q^{\dinv_m B} t^{\area_m(\De_{n,m}\setm B)}.
\eeq

\bfi
\begin{tikzpicture}
\draw (0,0)--(1,0) (1,0) grid (4,2) (2,2) grid (4,4) (3,4) grid (4,6) (4,6)--(4,8);
\draw[very thick] (0,0)--(1,0)--(1,2)--(2,2)--(2,4)--(3,4)--(3,6)--(4,6)--(4,8);
\end{tikzpicture}
\hs{50pt}
\begin{tikzpicture}
\fill[lightgray] (1,0) rectangle (2,2)  (2,3) rectangle (3,4) (3,4) rectangle (4,6);
\draw (0,0)--(1,0) (1,0) grid (4,2) (2,2) grid (4,4) (3,4) grid (4,6) (4,6)--(4,8);
\draw[very thick] (0,0)--(2,0)--(2,3)--(3,3)--(3,4)--(4,4)--(4,8);
\end{tikzpicture}
\hs{50pt}
\begin{tikzpicture}
\draw (0,0)--(1,0) (1,0) grid (4,2) (2,2) grid (4,4) (3,4) grid (4,6) (4,6)--(4,8);
\draw[very thick] (0,0)--(2,0)--(2,3)--(3,3)--(3,4)--(4,4)--(4,8);
\foreach \x in {0,.25,...,2} 
   \draw (\x,0) circle(.1);
\foreach \y in {0,.25,...,4}
   \draw (2,\y) circle(.1);
\foreach \x in {2,2.25,...,4} 
   \draw (\x,4) circle(.1);
\foreach \y in {4,4.25,...,8}
   \draw (4,\y) circle(.1);
\filldraw (2,2) circle(.1);
\filldraw (2,4) circle(.1);
\filldraw (4,6) circle(.1);
\filldraw (4,8) circle(.1);
\end{tikzpicture}
\capt{$\De_{4,2}$ and $(0,0,3,4)$  and its $2$-bounce path \label{Delta}}
\efi

The second definition of $C_{n,m}(q,t)$ which we will use involves the $m$-bounce statistic for boards $B\sbe\De_{n,m}$, so-called because it can be thought of as the path of a ball bouncing off of the boundary of $B$.  To define this statistic, we must first define the \emph{$m$-bounce path} $p$ for $B$ as follows.  Start $p$ from $(0,0)$ and take  horizontal steps stopping just short of the  first  lattice point interior to $B$.  (\emph{Interior} means strictly southeast of $B$'s boundary.)  Let $h_0$ denote the number of steps taken.  In the right-hand diagram of Figure~\ref{Delta}, the $2$-bounce path is shown in circles and $h_0=2$.  Now take $v_0=h_0$ vertical steps.  At the $i$th stage of this procedure, one starts from the current end of $p$ and moves horizontally $h_i$ steps, stopping just short of the first lattice point interior to $B$.  This is followed by $v_i$ vertical steps where
$$
v_i=h_i+h_{i-1}+\dots+h_{i-m+1}
$$
with the convention that $h_j=0$ if $j<0$.  In Figure~\ref{Delta} we mark the end of each vertical segment by a solid circle to make it clearer when there are horizontal segments of length zero.
One can prove that eventually this path will reach $(n,nm)$ where it terminates.  The complete sequence in Figure~\ref{Delta} is
$$
h_0 = 2,\quad v_0 = 2,\quad  h_1=0,\quad  v_1 = 0+2=2,\quad  h_2=2,\quad v_2=2+0=2,\quad
h_3=0,\quad v_3=0+2=2.
$$
Using this information, the \emph{$m$-bounce statistic} of $B\in\cB_{n,m}$ is
$$
\bounce_m B = \sum_{i\ge0} i h_i.
$$
So, returning to our example, 
$$
\bounce_2(0,0,3,4)=0\cdot 2 + 1\cdot 0+2\cdot 2 = 4.
$$
Using this statistic, we can define
\beq
\label{cB}
C_{n,m}(q,t)=\sum_{B\in\cB_{n,m}} q^{\area_m (\De_{n,m}\setm B)}t^{\bounce_m B}.
\eeq

To show that~\ree{cO} and~\ree{cB} define the same polynomials, we will use a bijection $\Phi:\cB_{n,m}\ra\cB_{n,m}$ introduced by Loehr~\cite{loe:csh}.  Given $B\in\cB_{n,m}$, let $n(\om_m(B))=(n_0,n_1,\dots)$.  To define $B'=\Phi(B)$, we will first define the bounce path $p$ of $B'$.   Starting with $\De_{n,m}$, draw the unique bounce path $p$ which has $h_i=n_i$ for all $i$.  It can be proved that such a path exists.  To illustrate, if  $\om=\om_2(0,0,3,4)=(0,2,1,2)$ then $n(\om)=(1,1,2,0)$.  So we would construct the bounce path with
$$
h_0=1,\quad v_0=1,\quad h_1=1,\quad v_1=1+1=2,\quad h_2=2,\quad v_2=2+1=3,\quad h_3=0,\quad v_3=0+2=2.
$$
This path is illustrated in Figure~\ref{Phi}.  Note that since $h_i=n_i$ for all $i$ we will have
$$
\bounce_m B' =\sum_{i \ge0} i h_i=\sum_{i\ge 0} i n_i = \area_m(\De_{n,m}\setm B)
$$
once $B'=\Phi(B)$ is defined.  Thus $\Phi$ will send the $\area_m$ statistic to the $\bounce_m$ statistic as desired to show the equality of the two definitions.

\bfi
\begin{tikzpicture}
\fill[lightgray] (1,0) rectangle (2,1);
\fill[lightgray] (2,1) rectangle (4,3);
\draw (0,0)--(1,0) (1,0) grid (4,2) (2,2) grid (4,4) (3,4) grid (4,6) (4,6)--(4,8);
\foreach \x in {0,.25,...,1} 
   \draw (\x,0) circle(.1);
\foreach \y in {0,.25,...,1}
   \draw (1,\y) circle(.1);
\foreach \x in {1,1.25,...,2} 
   \draw (\x,1) circle(.1);
\foreach \y in {1,1.25,...,3}
   \draw (2,\y) circle(.1);
\foreach \x in {2,2.25,...,4} 
   \draw (\x,3) circle(.1);
\foreach \y in {3,3.25,...,8}
   \draw (4,\y) circle(.1);
\filldraw (1,1) circle(.1);
\filldraw (2,3) circle(.1);
\filldraw (4,6) circle(.1);
\filldraw (4,8) circle(.1);
\end{tikzpicture}
\hs{100pt}
\begin{tikzpicture}
\draw (0,0)--(1,0) (1,0) grid (4,2) (2,2) grid (4,4) (3,4) grid (4,6) (4,6)--(4,8);
\draw[very thick] (0,0)--(1,0)--(1,1)--(2,1)--(2,2)--(3,2)--(3,3)--(4,3)--(4,8);
\draw (.3,.3) node {$0$} (.7,.7) node {$0$} (1.3,1.3) node {$1$} (1.7,1.7)
 node {$0$} (2.3,2.3) node {$2$} (2.7,2.7) node {$1$}  (3.3,3.3) node {$2$} (3.7,3.7) node {$2$} 
(3.7,4.5) node {$1$} (3.7,5.5) node {$2$} (3.7,6.5) node {$2$} (3.7,7.5) node {$2$}
;
\end{tikzpicture}
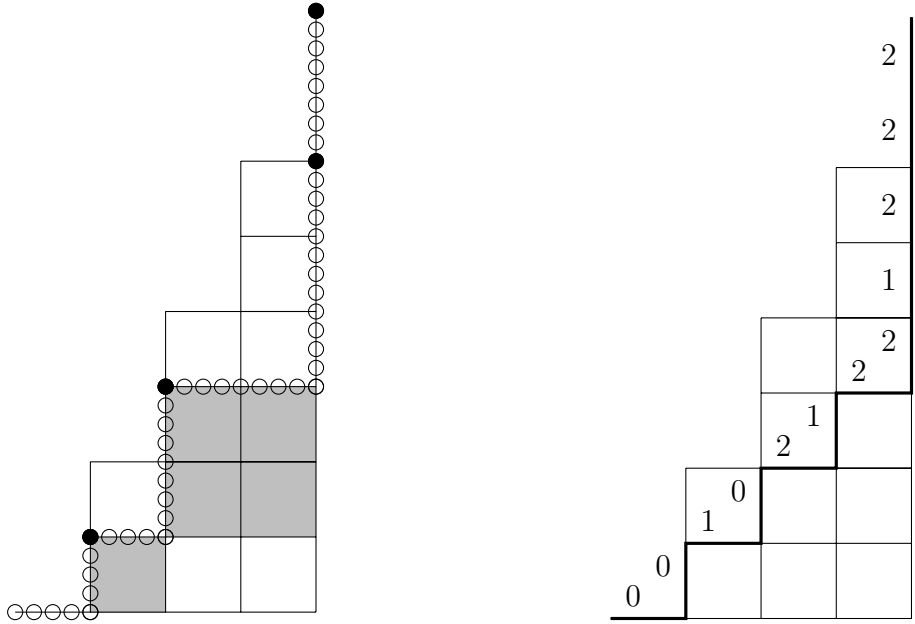
\capt{Constructing $\Phi(0,0,3,4)$ when $m=2$  \label{Phi}}
\efi

To obtain $B'$, we need one more concept.  Let $H_i$ be the set of steps counted by $h_i$ and similarly for $V_i$.   
It is convenient to let $V_{-1}$ and $H_{s+1}$ (where $V_s$ is the last vertical segment) consist of the vertices $(0,0)$ and $(n,nm)$, respectively. 
For $i\ge0$, define the \emph{$i$th bounce rectangle} $R_i$ of the bounce path $p$ to be the lattice rectangle whose southwest and northeast vertices are the first lattice point on $V_{i-1}$ and the last lattice point on $H_i$, respectively.  Note that if either of these sets of steps are empty, then the rectangle degenerates to a line.    On the left in Figure~\ref{Phi} the rectangles are $R_0$ which is the line segment from $(0,0)$ to $(1,0)$, $R_1$ which is a single square (shaded),  $R_2$ which is a $2\times 2$ rectangle (shaded), $R_3$ which is the line segment from $(4,3)$ to $(4,6)$, and $R_4$ which is the line segment from $(4,6)$ to $(4,8)$.  Now $B'$ will have bounce path $p$ if and only if the boundary of $B'$ travels  from the southwest corner to the northeast corner of each $R_i$ and coincides with the lowest step of $V_{i-1}$ on $p$.  (If $V_{i-1}=\emp$ then $R_i$ is a line segment and the bounce path just travels along it.)  To determine the part of the boundary of $B'$ in $R_i$, read the subword $\om^i$ of $\om$ consisting of the symbols $i,i-1,\dots,i-m$ from left to right, replacing each $i$ by an $E$ step and any of the other symbols by $N$ steps.  Continuing our example with $\om=(0,2,1,2)$, we have
$$
\om^0 = (0),\quad \om^1 = (0,1), \quad \om^2=(0,2,1,2),\quad \om^3=(2,1,2),\quad \om^4=(2,2)
$$
which translate to
$$
(E),\quad (N,E),\quad (N,E,N,E),\quad (N,N,N), \quad (N,N,N).
$$
The resulting path is displayed on the right in Figure~\ref{Phi} using thickened lines and labeled with the symbols from each $\om^i$.  The final result is the  board $B'=(0,1,2,3)=\Phi(0,0,3,4)$.  One can prove the following result.
\bth[\cite{loe:csh}]
\label{loe}
The map $\Phi:\cB_{n,m}\ra\cB_{n,m}$ is a bijection such that, if $\Phi(B)=B'$,

\vs{10pt}

\eqqed{
area_m(\De_{n,m}\setm B) =\bounce_m B' \qmq{and} \dinv_m B = \area_m (\De_{n,m}\setm B').
}
\eth

We now have everything in place to give the second proof of Theorem~\ref{card2}.

\bprf \emph{(Theorem~\ref{card2})}  Consider $\De_{N+1,m}$ where the first subscript has been chosen so that the whole $m$-weight equivalence class $[B]_m$ of $B$ is contained in $\cB_{N+1,m}$.  By Proposition~\ref{om}, we know that  $|[B]_m|$ is the number of rearrangements of $\om=\om_m(B)$ which give rise to a Ferrers board.  Note that rearranging $\om$ does not change the vector $n(\om)$.   So by the previous discussion, as $\om$ ranges over these rearrangements, $\Phi(B)$ ranges over all 
$B'\in \cB_{N+1,m}$ whose bounce path has horizontal components $h_i=n_i$ for all $i$. 

Now we can count all possible $B'$, rectangle by rectangle.  In $R_i$, the boundary of $B'$ must start with the first vertical step of $V_{i-1}$ and then travel to the last vertex of $H_i$.  To do this, it takes $h_i$ horizontal steps and 
$v_{i-1}-1=h_{i-1}+\dots+h_{i-m}-1$ vertical steps.  Remembering that $h_i=n_i$ for all $i$, we see that the number of possible boundary paths in $R_i$ is
$$
\binom{h_i+v_{i-1}-1}{h_i} = \binom{n_i+n_{i-1}+\dots+n_{i-m}-1}{n_i}
$$
Taking the product over $i$ gives the final result.
\eprf

%%%%%%%%%%%%%%%%%%%%%%%%%%%%%%%%
%
% 	Bounding boards
%
%%%%%%%%%%%%%%%%%%%%%%%%%%%%%%%%

\subsection{Bounding boards}
\label{bb}

Let $B=(b_1,\dots,b_n)$ be a Ferrers board with $b_1>0$.  We say that $\De_{N,m}$ is a \emph{bounding board} for $B$, or that $B$ \emph{fits} in $\De_{N,m}$ if $B'\sbe\De_{N,m}$ where $B'$ is $B$ padded with zeros until it has $N$ columns.  Usually we have been taking $N=|B|+1$ so that all members of $B$'s $m$-weight equivalence class will also fit in $\De_{N,m}$.  However, we can strengthen some of these results and make computations with examples easier if we relax this restriction on $N$.  For this, we need the following lemma.
\ble
\label{fit}
If $B$ fits in $\De_{N,m}$ for some $N$ and $B'\approx_m B$ then $B'$ fits in $\De_{N,m}$.
\ele
\bprf
It clearly suffices to prove this where $N$ is the smallest value such that $\De_{N,m}$ is a bounding board for $B$.  To find this value, we know that $B$ fits in $\De=\De_{|B|+1,m}$.  So  pad $B$ with zeros to get $B=(b_0,\dots,b_{|B|})$ and let $\om=\om_m(B)=(a_0,\dots,a_{|B|})$.    Let $c$ be the largest value such that $b_0=\dots=b_c=0$ and define the \emph{initial} and \emph{final} parts of $\om$ to be $\om_I=(a_0,\dots,a_c)$ and $\om_F=(a_{c+1},\dots,a_{|B|})$, respectively.  By definition of $c$, we have $\om_I=(0,m,\dots,cm)$ and 
\beq
\label{ac}
a_{c+1}<(c+1)m.
\eeq  
From the geometry of the situation we see that the smallest value of $N$ is determined by the minimum element in $\om_F$ and is, given that every change in $N$ results in an $m$-fold change in height,
$$
N=|B|+1-\lf a/m \rf_1.
$$
where $a$ is a minimum element in $\om_F$.

To complete the proof, we must show that if we make the same computations for $B'$, then we will get the same minimum value of $N$.  From the previous paragraph we see it suffices to show that $\om_F$ and $\om_F'$ have the same minimum value. 
There is no loss of generality in assuming $c\le c'$ which implies, since $\om'$ is a rearrangement of $\om$, that $\om_I\sbe\om_I'$ and $\om_F\spe\om_F'$ as sets.  Thus $a\le a'$ where $a=\min\om_F$ and $a'=\min\om_F'$.  Now suppose, towards a contradiction, that $a<a'$.  Since $a\not\in \om_F'$ we must have $a\in\om_I'$ and so $a=im$ for some $i$.  If $i\le c$ then $a$ appears in both $\om_I$ and $\om_F$.  And since $a=im$ appears only once in $\om_I'$, it must also appear in $\om_F'$ which contradicts the fact that the $a'=\min\om_F'$.  If $i>c$ then $a=im\ge(c+1)m$.  But then, by equation~\ree{ac}, $a_{c+1}<a$ which contradicts the fact that $a=\min\om_F$.  So in either case we have a contradiction and are done with the proof of the lemma.
\eprf

From now on we make the convention that if we say that $B$ fits in $\De_{N,m}$ and then compute $\om_B(B)$ we are using the version of $B$ with $N$ columns.  Even though the previous lemma makes it possible to choose different values of $N$, it is important to be consistent and keep the same value during a given argument.  We illustrate the use of the lemma with the following result.   An example follows the proof.
\bpr
The set of boards fitting in  $\cB_{n,m}$ is a union of $m$-weight equivalence classes.  The number of such classes is $(m+1)^{n-1}$.
\epr
\bprf
The first statement follows immediately from Lemma~\ref{fit}.  For the second,  applying $\Phi$  as in the second proof of Theorem~\ref{card2}, we see that this amounts to counting all possible bounce paths $p$ in $\cB_{n,m}$.  Note that such a $p$ is completely determined by the composition $h=(h_0,h_1,\dots,h_s)$ of $n$ where $h_s$ is the last positive $h_i$.  (Any $h_i=0$ for $i>s$ do not affect the path since it will have reached the line $x=n$ and so must go vertically up to $(n,mn)$ no matter how many more zero components there are.)  We also have $h_0>0$ since $\De_{n,m}$ starts with a column of zero height, forcing $p$ to initially move at least one  unit horizontally.  Another fact about $h$ is that it can contain at most $m-1$ consecutive parts equal to zero.  (See~\cite{loe:csh} for a proof.)   These are the only restrictions on $h$ and every $h$ satisfying them comes from a bounce path.

We count the desired compositions using the usual ``slashes and dashes" method.  Consider a string of $n$ ones which determine $n-1$ spaces in between each pair of consecutive ones.  In each space we put one of the following: nothing, a plus sign, or up to $m-1$ zeros.  This will give a composition by adding together the ones with plus signs in between and considering them and the inserted zeros as the parts of the composition.  So in each of the $n-1$ spaces we have $m+1$ choices of what to insert, giving a total count of $(m+1)^{n-1}$.
\eprf

To illustrate the proof, consider the case when $m=3$, $n=8$, and $h=(2,3,0,2,0,0,1)$.  This $h$ would be obtained from a string of $8$ ones by inserting elements to form 
$$
1+1\ \ 1+1+1\ \ 0\ \ 1+1\ \ 0\ \ 0\ \ 1.
$$

%%%%%%%%%%%%%%%%%%%%%%%%%%%%%%%%
%
% 	A $q$-analogue
%
%%%%%%%%%%%%%%%%%%%%%%%%%%%%%%%%

\subsection{A $q$-analogue}
\label{qa}

We next present a $q$-analogue of Thoerem~\ref{card2}.  For its statement, we need the \emph{$q$-binomial coefficients} defined by
$$
\gauss{n}{k}_q=\frac{[n]_q!}{[k]_q![n-k]_q!}
$$
where $[n]_q!$ is the result of setting $p=1$ in $[n]_{p,q}!$.
\bth
\label{card2q}
Let $B_0$ be a Ferrers board fitting inside $\De_{N,m}$ and let
$n(\omega_m(B_0))=(n_0,n_1,\ldots)$. We have
$$ 
\sum_{ B\approx_m B_0} q^{\dinv_m B} =
q^{c(n_0,n_1,\ldots)}  \prod_{i\geq 1}\gauss{n_i+n_{i-1}+\cdots+n_{i-m}-1}{n_i}_q,
$$
where 
$$
c(n_0,n_1,\ldots)=m\sum_{i\geq 0} \binom{n_i}{2}
 +\sum_{i\geq 1} n_i\sum_{j=1}^m (m-j)n_{i-j}.
$$
\eth
\bprf  Using Theorem~\ref{loe} and arguing as in the first paragraph of the second proof of Theorem~\ref{card2}, we see that the summation side of the identity in the theorem is
$$
\sum_{B'} q^{\area_m(\De_{N,m}\setm B')}.
$$
where the sum is over all $B'$ that have the  bounce path $p$ with $h_i=n_i$ for all $i$.  Let $a$ be the number of squares between $p$ and $\De_{N,m}$ and $a(B')$ be the number of squares between $B'$ and $p$ so that 
$\area_m(\De_{N,m}\setm B')=a+a(B')$.  The theorem will follow from the following two claims.

First, we claim that 
$$
q^a=q^{c(n_0,n_1,\ldots)}
$$
It suffices to show that the number of squares  below $\De_{N+1,m}$ which are above the horizontal steps of $p$ corresponding to  $n_i$ is $m\binom{n_i}{2} +  n_i\sum_{j=1}^m (m-j)n_{i-j}$.  This area is broken into a triangluar region above the line 
$y=m(n_{i-1}+\dots+n_{i-m})x$ and a rectangular region below.  Easy calculations now show that $m\binom{n_i}{2}$ is the number of cells in the triangular portion while $n_i\sum_{j=1}^m (m-j)n_{i-j}$ gives the analogous count for the rectangle.  This proves the first claim.

The second claim is that
$$
\sum_{B'} q^{a(B')}= \gauss{n_i+n_{i-1}+\cdots+n_{i-m}-1}{n_i}.
$$
For this it suffices to show that the $i$th factor is the generating function for the possible areas in the $i$th rectangle, $R_i$.  Recall that the boundary of $B'$ in this rectangle starts with a vertical step and then ranges over all lattice paths with $n_i$ horizontal steps and $n_{i-1}+n_{i-2}+\dots+n_{i-m}-1$ vertical steps.  It is well known that this count is given by the $q$-binomial coefficient above.  So we are done with the second claim and the second proof of the theorem.
\eprf

Since the $q$-binomial coefficients have leading and constant coefficients equal to one, we immediately get two distinguished representatives of an $m$-weight equivalence class from the previous theorem.
\bco
\label{B_1B_2}
Every $m$-weight equivalence class contains unique boards $B_1$ and $B_2$
such that
$$
\dinv_m B_1\leq\dinv_m B\leq \dinv_m B_2
$$
for all boards $B$ in the equivalence class.  
\eco

Note that, from the proof of Theorem~\ref{card2q}, $B_1$ and $B_2$ are uniquely defined by the fact that $B_1'=\Phi(B_1)$ and $B_2'=\Phi(B_2)$ are the most northwest and most southeast possible boards with the fixed bounce path, $p$ for their equivalence class.   So the boundary of $B_1'$ coincides with $p$.  And the boundary of $B_2'$ coincides with the southern most step on each vertical portion of $p$ and then proceeds to the next such step by taking a sequence of horizontal steps followed by a sequence of vertical steps.  

There is a simpler way to compute $B_1$ and $B_2$ which we now describe.  Suppose $B$ is a board in the $m$-weight equivalence class in question and suppose that $\De_{N,m}$ is any bounding board for $B$.  Let $\om=\om_m(B)$ and let $\omh$ be the unique weakly increasing rearrangement of $\om$.  In other words, $\omh$ is the rearrangement considered in the proof of Theorem~\ref{mres}.  In that demonstration we proved that $\omh=\om(\Bb_1)$ for some Ferrers board $\Bb_1$.    Continuing our running example with $B=(0,0,3,4)$ we have $\om=\om_2(B)=(0,2,1,2)$ and so $\omh=(0,1,2,2)$.  This corresponds to the board $\Bb_1=(0,1,2,4)$.  The reader can now verify that $\Phi(\Bb_1) = (0,1,3,3)$ whose boundary coincides with the bounce path  of $\Phi(B)$.  In general, it is easy to see from the definition of $\Phi$ that this last statement holds for any $B$ and so $\Bb_1=B_1$.

We now want to indicate another construction of $\Bb_1$ which will more closely parallel the construction of the board $\Bb_2$ which will equal $B_2$.  Let $n=n(\om)=(n_0,n_1,\dots,n_t)$.  We will construct a sequence of vectors $\omh_0,\dots,\omh_t$ as follows.  Let $\omh_0$ consist of $n_0$ zeros.  Once $\omh_{i-1}$ has been created, we obtain $\omh_i$ by inserting $n_i$ copies of $i$ in the position as far to the right in $\omh_{i-1}$ such that $\omh_i$ still satisfies the conditions of Proposition~\ref{oequiv} and so corresponds to a board.  It is easy to see by induction on $i$ that this amounts to putting the $i$'s at the end of $\omh_{i-1}$ and so $\omh_t=\omh$ as defined in the previous paragraph.  In the running example
$$
\omh_0=(0),\ \omh_1=(0,1),\ \omh_2=(0,1,2,2)=\omh.
$$
Now to construct $\Bb_2$, we build a sequence $\omc_0,\dots,\omc_t$ where $\omc_0=\omh_0$  and $\omc_i$ is obtained from $\omc_{i-1}$  by inserting the $n_i$ copies of $i$ as far to the left as possible so that $\omc_i$ still satisfies Proposition~\ref{oequiv}.  Finally, we let $\Bb_2$ be the board with $\om_m(\Bb_2)=\omc_t$.  
Using our example once more, we have
$$
\omc_0=(0),\ \omc_1=(0,1),\ \omc_2=(0,2,2,1)
$$
so that $\Bb_2=(0,0,2,5)$.  Again, the reader can verify that $\Phi(\Bb_2)=(0,1,2,2)$ whose boundary is the most southeastern having the same bounce path as $\Phi(B)$ both in this example and in general.

We will now use the descriptions of $\Bb_1$ and $\Bb_2$ to give a second proof of Corollary~\ref{B_1B_2} without using the map $\Phi$.

\bprf \emph{(Corollary~\ref{B_1B_2})}
We show that $\dinv_m \Bb_1$ in the unique minimum in the $m$-weight equivalence class of $B$.  The proof for $\Bb_2$ is similar.  It suffices to show that $\omh$ as defined above has $\dinv_m\omh$ as a unique minimum among all rearrangements of $\om$.  Let  $\pi=a_1\dots a_n$ be any sequence of real numbers with $a_i>a_{i+1}$ and let $\pib$ be $\pi$ with $a_i$ and $a_{i+1}$ interchanged.  We will show that $\dinv_m\pib\le \dinv_m\pi$, and that under certain circumstances $\dinv_m\pib< \dinv_m\pi$.  Finally we will demonstrate that we can use adjacent transpositions to transform any rearrangement of $\om$ into $\omh$ while weakly decreasing $\dinv_m$ at each step and strictly decreasing the statistic during at least one step.  Note that it does not matter whether the intermediate rearrangements of $\om$ correspond to Ferrers boards since we only need to prove the inequalities for $\dinv_m$.

By considering the contribution of each pair $a_i,a_j$ in $\pi$, one can rewrite the definition of  $\dinv_m\pi$ in equation~\ree{dinvdef} as
$$
\dinv_m \pi =\sum_{i<j} f_m(a_i-a_j)
$$
where
$$
f_m(d)=
\begin{cases}
m-d+1	&\text{if $0< d\le m$,}\\
m+d     &\text{if $-m\le d\le 0$,}\\
0	&\text{otherwise.}
\end{cases}
$$
Considering the $\pi$ and $\pib$ defined above, we see that all terms of $\dinv_m\pi-\dinv_m\pib$  cancel except for those corresponding to the given value of $i$ and $j=i+1$.   Letting $d=a_i-a_{i+1}>0$ we see that if $d>m$ then $\dinv_m\pi-\dinv_m\pib=0$, and  if 
$0<d\le m$ then
$$
\dinv_m\pi-\dinv_m\pib = f_m(d)-f_m(-d)=(m-d+1)-(m+(-d))=1.
$$
So in either case  $\dinv_m\pib\le \dinv_m\pi$, and if $0<d\le m$ then $\dinv_m\pib< \dinv_m\pi$.  Now given any rearrangement $\om'$ of $\om$, we can transform $\om'$ into $\omh$ by eliminating adjacent inversions at each step and thus always weakly decreasing $\dinv_m$.  Furthermore, by condition (ii) of Proposition~\ref{oequiv}, the last step will strictly decrease $\dinv_m$ because $\omh$ corresponds to a Ferrers board.  This completes the proof.
\eprf

As a final remark about this section, the reader will have noticed that we gave two proofs of Theorem~\ref{card2}, one from first principles and one using the $q,t$-Catalan machinery.  However, we only presented a proof of the latter type for its $q$-analogue, Theorem~\ref{card2q}.   This is because  a demonstration of the former type already exists in the literature, see the proof of equation~(13) in~\cite{loe:csh} and set $t=1$.

%%%%%%%%%%%%%%%%%%%%%%%%%%%%%%%%
%
% 	OPEN QUESTIONS
%
%%%%%%%%%%%%%%%%%%%%%%%%%%%%%%%%

\section{Open questions}
\label{oq}

\subsection{Counting $m$-level equivalence classes}

One would like to use Theorem~\ref{mft} to obtain a  formula for the size of any $m$-level equivalence class.  But the extra term in the factor for columns at the end of their zone may make this difficult to do.   Perhaps one could at least find a formula by putting some extra condition on the remainder of the zones of the board such as we have done by imposing the singleton restriction.

\subsection{A weighted $p,q$-analogue}

It would be very interesting to find a $p,q$-analogue of the $m$-weight Factorization Theorem, Theorem~\ref{mwft}.  The main stumbling block seems to be finding the correct way to translate equation~\ree{y} which is necessary since the final result cannot depend on the values of the $y_i$.  In particular, one would need a way of writing $[n-1]$ as a linear combination of the $[y_i]$ where the coefficient of $[y_i]$ did not depend on the other $[y_j]$.  It is not clear how to do this.

\subsection{$p,q$-hit and $q$-hit numbers for Ferrers boards}

Given a board $B$ contained in the $n \times n$ board, 
the \emph{$k$-th hit number} of $B$ with respect to $n$, $h_{k,n}(B)$, is 
defined be the number of permutations $\sigma \in \fS_n$ such that 
the rook placement corresponding to $\sigma$ has exactly $k$ rooks 
in $B$.  A classical result of 
Riordan and Kaplansky \cite{kr:pra} gives a simple 
relationship between the hit numbers and the rook numbers of 
$B$, namely, 
$$\sum_{k=0}^n h_{k,n}(B)x^k = \sum_{k=0}^n r_k(B)(n-k)! (x-1)^k.$$
There is a natural analogue of the hit numbers which we call 
the $m$-level hit numbers for 
boards contained in the $mn \times n$ board. That is, 
Briggs and Remmel~\cite{br:mrn} observed that a 
$m$-level rook placement of $n$ rooks in the 
$mn \times n$ board can naturally be identified with 
an element in the wreath product $C_m \wr \fS_n$ of the cyclic 
group $C_m$ with the symmetric group $\fS_n$. 
Thus given a board $B$ contained in the $mn \times n$ board, we can 
define the \emph{$k$-th $m$-level hit number} of $B$ with respect to $n$, 
$h_{k,n,m}(B)$, to be the number of $\sigma \in C_m \wr \fS_n$ such that 
the rook placement corresponding to $\sigma$ has exactly $k$ rooks 
in $B$. There is a natural extension of the Riordan 
and Kaplansky's result which relates the $m$-level hit numbers 
of $B$ to the $m$-level rook numbers of $B$, namely, 
$$
\sum_{k=0}^n h_{k,n,m}(B)x^k = \sum_{k=0}^n r_{k,m}(B) (m(n-k))\downarrow_{n-k,m} (x-1)^k.
$$

Briggs and Remmel~\cite{br:mrn} defined a $p,q$-analogue of 
the $m$-level hit numbers $h_{k,n,m}[B]$ by the following equation: 
\begin{equation}\label{pqhit}
\sum_{k=0}^n h_{k,n,m}[B]x^k = \sum_{k=0}^n r_{k,m}[B] [m(n-k)]\downarrow_{n-k,m} p^{m(\binom{k+1}{2} + k(n-k))} \prod_{\ell = n-k+1}^n 
(x-q^{m\ell}p^{m(n-\ell)}).
\end{equation}
They proved that for all singleton boards $B$, $h_{k,m,n}[B]$ is a polynomial 
in $p$ and $q$ with non-negative integer coefficients. It is 
natural to ask whether such a result can be extended to all 
Ferrers boards. The answer is no, in general. For example, suppose $m=2$ and 
$B = (1,1,1)$.  Then it is easy to check that $r_{0,2}[B] = q^3$,  
$r_{1,2}[B] = p^{-2} +p^{-4}q+p^{-6}q^2$, and 
$r_{2,2}[B] =r_{3,2}[B] =0$. Thus in this case, (\ref{pqhit}) becomes 
\begin{eqnarray*}
\sum_{k=0}^3 h_{k,3,2}[B]x^k &=& 
r_{0,2}[B] [6][4][2]+r_{1,2}[B][4][2]p^6(x-q^6) \\
&=&[4][2](q^3[6]+(p^4+p^2q+q^2)(x-q^6))
\end{eqnarray*} 
so that
\begin{eqnarray*}
 h_{0,3,2}[B] &=& [4][2](q^3[6]-(p^4q^6+p^2q^7 +q^8)), \\
h_{1,3,2}[B] &=& [4][2] (p^4+qp^2 +q^2), \ \mbox{and} \\
h_{2,3,2}[B] &=& h_{3,2,3}[B] =0.
\end{eqnarray*}
Note that $q^3[6]-(p^4q^6+p^2q^7 +q^8)$ does not have non-negative 
coefficients since the terms in $q^3[6]$ are all homogeneous of degree 
$8$ and, hence, there is nothing to cancel $-p^4q^6$ or $-p^2q^7$. In 
However, if we set $p=1$  then we do end up with 
$h_{k,3,2}[B]$ being a polynomial in 
$q$ with non-negative integer coefficients for all $k$.  Thus a natural 
question to ask is whether there are non-singleton 
Ferrers boards $B$ such that $m$-level hit numbers 
$h_{k,n,m}[B]$ are always polynomials in $p$ and $q$ with 
non-negative integer coefficients and, if so, can one classify 
such non-singleton Ferrers boards.  Similarly, it would be 
interesting to classify those Ferrers boards $B$ such that 
$h_{k,n,m}[B]$ are polynomials in $q$ with non-negative coefficients
when $p=1$.

\end{document}